\documentclass[reqno,11pt,letterpaper]{amsart}
\usepackage{mathrsfs}
\usepackage{amssymb}

\usepackage{amsthm}
\usepackage{amsmath}
\usepackage{amsfonts}
\usepackage{bm, enumerate,xcolor}

\usepackage{graphicx}

\numberwithin{equation}{section}

\newtheorem{theorem}{Theorem}[section]
\newtheorem{lemma}[theorem]{Lemma}
\newtheorem{corollary}[theorem]{Corollary}
\newtheorem{proposition}[theorem]{Proposition}

\theoremstyle{definition}
\newtheorem{definition}[theorem]{Definition}
\newtheorem{assumption}[theorem]{Assumption}
\newtheorem{example}[theorem]{Example}

\theoremstyle{remark}
\newtheorem{remark}[theorem]{Remark}

\allowdisplaybreaks[1]

\begin{document}

\title[Steady compressible self-similar flows]{The steady inviscid compressible self-similar flows and the stability analysis}

\author{Shangkun WENG}
\address{School of mathematics and statistics, Wuhan University, Wuhan, Hubei Province, China, 430072.}
\email{skweng@whu.edu.cn}


\author{Hongwei YUAN}
\address{Department of Mathematics, University of Macau, Taipa, Macau, China.}
\email{hwyuan@um.edu.mo}

\date{}
\maketitle

\def\be{\begin{eqnarray}}
\def\ee{\end{eqnarray}}
\def\ba{\begin{aligned}}
\def\ea{\end{aligned}}
\def\bay{\begin{array}}
\def\eay{\end{array}}
\def\bca{\begin{cases}}
\def\eca{\end{cases}}
\def\p{\partial}
\def\hphi{\hat{\phi}}
\def\bphi{\overline{\phi}}
\def\no{\nonumber}
\def\eps{\epsilon}
\def\de{\delta}
\def\De{\Delta}
\def\om{\omega}
\def\Om{\Omega}
\def\f{\frac}
\def\th{\theta}
\def\vth{\vartheta}
\def\la{\lambda}
\def\lab{\label}
\def\b{\bigg}
\def\var{\varphi}
\def\na{\nabla}
\def\ka{\kappa}
\def\al{\alpha}
\def\La{\Lambda}
\def\ga{\gamma}
\def\Ga{\Gamma}
\def\ti{\widetilde}
\def\wti{\widetilde}
\def\wh{\widehat}
\def\ol{\overline}
\def\ul{\underline}
\def\Th{\Theta}
\def\si{\sigma}
\def\Si{\Sigma}
\def\oo{\infty}
\def\q{\quad}
\def\z{\zeta}
\def\co{\coloneqq}
\def\eqq{\eqqcolon}
\def\di{\displaystyle}
\def\bt{\begin{theorem}}
\def\et{\end{theorem}}
\def\bc{\begin{corollary}}
\def\ec{\end{corollary}}
\def\bl{\begin{lemma}}
\def\el{\end{lemma}}
\def\bp{\begin{proposition}}
\def\ep{\end{proposition}}
\def\br{\begin{remark}}
\def\er{\end{remark}}
\def\bd{\begin{definition}}
\def\ed{\end{definition}}
\def\bpf{\begin{proof}}
\def\epf{\end{proof}}
\def\bex{\begin{example}}
\def\eex{\end{example}}
\def\bq{\begin{question}}
\def\eq{\end{question}}
\def\bas{\begin{assumption}}
\def\eas{\end{assumption}}
\def\ber{\begin{exercise}}
\def\eer{\end{exercise}}
\def\mb{\mathbb}
\def\mbR{\mb{R}}
\def\mbZ{\mb{Z}}
\def\mc{\mathcal}
\def\mcS{\mc{S}}
\def\ms{\mathscr}
\def\lan{\langle}
\def\ran{\rangle}
\def\lb{\llbracket}
\def\rb{\rrbracket}
\def\fr#1#2{{\frac{#1}{#2}}}
\def\dfr#1#2{{\dfrac{#1}{#2}}}
\def\u{{\textbf u}}
\def\v{{\textbf v}}
\def\w{{\textbf w}}
\def\d{{\textbf d}}
\def\nn{{\textbf n}}
\def\x{{\textbf x}}
\def\e{{\textbf e}}
\def\D{{\textbf D}}
\def\U{{\textbf U}}
\def\M{{\textbf M}}
\def\F{{\mathcal F}}
\def\I{{\mathcal I}}
\def\W{{\mathcal W}}
\def\div{{\rm div\,}}
\def\curl{{\rm curl\,}}
\def\R{{\mathbb R}}
\def\FF{{\textbf F}}
\def\A{{\textbf A}}
\def\R{{\textbf R}}
\def\r{{\textbf r}}

\begin{abstract}
  We investigate the steady inviscid compressible self-similar flows which depends only on the polar angle in spherical coordinates. It is shown that besides the purely supersonic and subsonic self-similar flows, there exists purely sonic flows, Beltrami flows with a nonconstant proportionnality factor and smooth transonic self-similar flows with large vorticity. For a constant supersonic incoming flow past an infinitely long circular cone, a conic shock attached to the tip of the cone will form, provided the opening angle of the cone is less than a critical value. We introduce the shock polar for the radial and polar components of the velocity and show that there exists a monotonicity relation between the shock angle and the radial velocity, which seems to be new and not been observed before. If a supersonic incoming flow is self-similar with nonzero azimuthal velocity, a conic shock also form attached to the tip of the cone. The state at the downstream may change smoothly from supersonic to subsonic, thus the shock can be supersonic-supersonic, supersonic-subsonic and even supersonic-sonic where the shock front and the sonic front coincide. We further investigate the structural stability of smooth self-similar irrotational transonic flows and analyze the corresponding linear mixed type second order equation of Tricomi type. By exploring some key properties of the self-similar solutions, we find a multiplier and identify a class of admissible boundary conditions for the linearized mixed type second-order equation. We also prove the existence and uniqueness of a class of smooth transonic flows with nonzero vorticity which depends only on the polar and azimuthal angles in spherical coordinates.
\end{abstract}

\section{Introduction}

We focus on the study of three-dimensional steady compressible Euler flows past an infinitely long circular cone. These steady inviscid compressible flows are governed by the following system of equations:
\begin{align}\label{comeuler3d}
\begin{cases}
\div(\rho{\bf u})=0,\\
\div(\rho {\bf u}\otimes{\bf u})+\nabla p=0,\\
\div(\rho E{\bf u}+p{\bf u})=0,\\
\end{cases}
\end{align}
where ${\bf u}=(u_1,u_2,u_3)$ is the velocity, $\rho$ is the density. The pressure $p$ can be expressed as $p(\rho,S)=A(S)\rho^{\gamma}$; here $A(S) = a e^{\frac{S}{c_v}}$ is a one-to-one correspondence between the entropy $S$ and the quantity $A$, given fixed positive constants $a$ and $c_v$ as well as the exponent $\gamma>1$; we also refer to the quantity $A$ as entropy, which will not create any ambiguity. The total energy is defined as $E=\frac{1}{2} |{\bf u}|^2 + e(\rho,S)$, where $e(\rho,S)=\frac{p}{(\gamma-1)\rho}$ represents the internal energy. Additionally, we denote the enthalpy and Bernoulli's function as
\begin{align}\label{Bernoulli}
h(\rho,S)=\frac{\gamma p(\rho,S)}{(\gamma-1)\rho}\text{ and }B=\frac{1}{2} |{\bf u}|^2+ h(\rho,S),
\end{align}
respectively.

To serve our purpose, we introduce the spherical coordinates $(r, \theta, \varphi)$. 
\begin{align}\no
x_1=r\sin\theta \cos\varphi,\ \ \  x_2=r\sin\theta \sin\varphi,\ \ \ \ x_3=r\cos\theta.
\end{align}
We can decompose the velocity vector $\mathbf{u}$ in spherical coordinates as follows:
\begin{equation}\no
\mathbf{u} = U_1 \mathbf{e}_r + U_2 \mathbf{e}_\theta + U_3 \mathbf{e}_\varphi,
\end{equation}
where $U_1$, $U_2$, and $U_3$ represent the radial, polar, and azimuthal components of the velocity, respectively. Here, $\mathbf{e}_r$, $\mathbf{e}_\theta$, and $\mathbf{e}_\varphi$ are the unit vectors in the radial, polar, and azimuthal directions, respectively
\begin{align}\no
{\bf e}_r=\begin{pmatrix} \sin\theta \cos\varphi\\ \sin\theta\sin\varphi \\\cos\theta \end{pmatrix},
{\bf e}_\theta=\begin{pmatrix} \cos\theta\cos\varphi\\ \cos\theta \sin\varphi\\-\sin\theta\end{pmatrix},
{\bf e}_\varphi=\begin{pmatrix}-\sin\varphi\\ \cos\varphi\\0\end{pmatrix}.
\end{align}
By expressing the velocity vector in this decomposed form, for $C^1$-smooth solutions, the system \eqref{comeuler3d} can be rewritten as
\begin{align}\label{SCEQF}
\begin{cases}
\p_r (\rho U_1)+ \frac1r\p_{\theta}(\rho U_2)+\frac1{r\sin\theta}\p_{\varphi}(\rho U_3)+\frac2r \rho U_1+\frac1{r\tan\theta}\rho U_2=0,\\
(U_1\p_r+\frac{U_2}{r}\p_{\theta}+ \frac{U_3}{r\sin\theta}\p_{\varphi}) U_1+\frac1\rho \p_r p-\frac{U_2^2+U_3^2}r=0,\\
(U_1\p_r+\frac{U_2}{r}\p_{\theta}+ \frac{U_3}{r\sin\theta}\p_{\varphi}) U_2+\frac1{r\rho} \p_{\theta}p+\frac{U_1U_2}r-\frac{U_3^2}{r\tan\theta}=0,\\
(U_1\p_r+\frac{U_2}{r}\p_{\theta}+ \frac{U_3}{r\sin\theta}\p_{\varphi}) U_3+\frac1{r\rho\sin\theta} \p_{\varphi}p+\frac{U_1U_3}r+\frac{U_2U_3}{r\tan\theta}=0,\\
(U_1\p_r+\frac{U_2}{r}\p_{\theta}+ \frac{U_3}{r\sin\theta}\p_{\varphi}) A=0.
\end{cases}
\end{align}

Consider the self-similar solutions to the system of equations \eqref{comeuler3d}, then all the fluid parameters depend only on the polar angle, we can express the velocity vector as ${\bf u}=U_1(\theta){\bf e}_r+U_2(\theta){\bf e}_\theta+U_3(\theta){\bf e}_\varphi$, the density $\rho=\rho(\theta)$ and the pressure $p=A(\theta)\rho^\gamma$. Under these assumptions, the system of equations \eqref{SCEQF} reduces to:
\begin{eqnarray}\label{SCEQF1}
\begin{cases}
(\rho U_2)'(\theta)+2 \rho U_1+\frac1{\tan\theta}\rho U_2=0,\\
U_2 U_1'-(U^2_2+U_3^2)=0,\\
U_2U_2'+\frac1\rho \frac{d}{d\theta}p(\rho)+U_1U_2-\frac{U_3^2}{\tan\theta}=0,\\
U_2U_3'+U_1U_3+\frac{U_2U_3}{\tan\theta}=0,\\
U_2 A'=0.
\end{cases}
\end{eqnarray}

The studies of the self-similar solution to \eqref{SCEQF} depending only on the polar angle has a long history, at least starting from 1920s \cite{bus31,ma37,tm33}. This kind of solutions had been studied by Taylor and Maccoll \cite{tm33} within the irrotational flows where $U_3\equiv 0$. For the constant supersonic incoming flow past an infinitely long circular cone, Busemann \cite{bus31} introduced the shock polar to describe the one-parametric family of the possible flow state at downstream after crossing the shock. The cylindrical coordinates were used in \cite{bus31,cf48} and an ODE system was derived to describe the evolution of the flow at downstream. If the shock is prescribed, the end state on the cone surface will form a curve which Busemann called it ``apple curve". If the opening angle of the cone is less than some critical value, from the apple curve there may exist two possible shock solutions: supersonic-supersonic weak shock and supersonic-subsonic strong shock. One may refer to the classical book by Courant and Friedrichs \cite[Sections 121, 155, 156]{cf48} for more detailed description. The question that which solution is more likely to occur in reality is a big challenging problem in gas dynamics and had draw favorable attentions of many mathematicians \cite{cl00,ch01,cxy02,el08,ll99,maj90}. Recently, using the spherical coordinates, Lien-Liu-Peng \cite{llp22} considered the problem the supersonic incoming flow past the cone by prescribing the initial data on the cone and integrating the ODEs until the Rankine-Hugoniot condition is satisfied on the shock cone. They further proved the existence of self-similar irrotational smooth transonic flows for suitably chosen incoming flows, even for the small opening angle case.

Here our consideration to \eqref{SCEQF1} will cover both the irrotational and rotational flows. We first consider the initial value problem to \eqref{SCEQF1}. It is found that besides the purely supersonic and subsonic flows, there exist purely sonic flows and Beltrami flows with a nonconstant proportionnality factor. There also exist smooth transonic self-similar flows with large vorticity, where the sonic points are noncharacteristic degenerate and form a conic surface. We further concern the shock solution to \eqref{SCEQF1} which corresponds to the conic shock formed when supersonic incoming flows past an infinitely long circular cone. Different from \cite{cf48,xy08} where they used the cylindrical coordinates, the Rankine-Hugoniot jump conditions in spherical coordinates take the same form as the normal shock jump conditions in a duct. We introduce a new shock polar for $U_1$ and $U_2$ and find that there exists a monotonicity relation between the shock angle and the radial velocity for the case of a constant supersonic incoming flow past a cone, which seems to be new and not been observed before. For a supersonic incoming flow with nonzero swirl velocity past an infinitely long circular cone, if the opening angle of the cone belong to some set, a conic shock attached to the tip of the cone will form. The flow state at downstream may change smoothly from supersonic to subsonic, thus the shock can be supersonic-supersonic, supersonic-subsonic and even supersonic-sonic where the shock front and the sonic front coincide. Similar phenomena happen for the radially symmetric spiral flows in an annulus, which had been studied by the authors in \cite{wxy21a}.

There have many literatures working on the stability analysis of the self-similar solutions to \eqref{SCEQF}. The potential flow model is widely used in many circumstances for the studies of supersonic flows \cite{cf48,go97,mt87}. The local existence of supersonic flows past a cone for the potential flow were proved in \cite{cl00} for the symmetric case and \cite{ch01} for the nonsymmetric case. Lien-Liu in \cite{ll99} constructed a global weak solution via a modified Glimme scheme and obtained the long time asymptotic behavior under some smallness assumptions of the sharp vertex angle and the shock strength. Based on some global uniform weighted energy estimates, Chen-Xin-Yin in \cite{cxy02} proved the stability of steady supersonic conic shock for steady supersonic flow past an infinite curved and symmetrically perturbed cone, provided the vertex angle of the cone is suitable small and the Mach number is sufficiently large. The authors in \cite{cy07,cy09} improved the global existence result for the symmetrically perturbed supersonic flow under the assumption that the vertex angle less than a critical value. Later on, Li-Ingo-Yin \cite{liy14} established the global existence of a three dimensional supersonic conic shock solution for the potential equation by exploring a new Hardy-type inequality on the shock surface. It should be emphasized that except the main curved conic shock, there are no other discontinuities in the solutions constructed in \cite{cxy02,cy07,cy09,liy14,xy06}, which is quite different from the results in \cite{ll99,wz09}. The authors in \cite{ckz21} established the structural stability of conical shocks in three dimensional steady supersonic isothermal potential flows with axisymmetry past Lipschitz perturbed cones.

For the stability of transonic conic shock past a cone, under the assumption that the velocity of incoming flow is appropriately large, the authors in \cite{xy08,xy09} established the stability of the self-similar transonic shock wave for the symmetrically perturbed supersonic flow past an infinitely long cone when the pressure at infinity is appropriately large. The transonic conic shock was shown to be unstable in \cite{xy10} for the steady Euler system due to the influence of the rotation. For supersonic flow past 2-D wedge or 3-D wing, one may refer to \cite{ccf16,cf17,czz06,ccx21,cf07,el08,lxy15,wz09,yin06} and the references therein for more details.

In this paper, we will investigate the structural stability of the smooth irrotational transonic self-similar flows under small perturbations of suitable boundary conditions. Wang and Xin \cite{wx19,wx21} established the existence and uniqueness of Meyer type irrotational transonic flows satisfying some physical boundary conditions on the De Laval nozzle by transforming the two dimensional steady irrotational compressible Euler equations to the Chaplygin equation on the potential function and stream function plane. The authors in \cite{wxy21b} had investigated the structural stability of a class of radially symmetric transonic flows with nonzero angular velocity under small perturbations of suitable boundary conditions and established the existence and uniqueness of smooth transonic spiral flows with nonzero vorticity in concentric cylinders based on the analysis of a linear mixed type second order equation of Tricomi type and the deformation-curl decomposition developed in \cite{wex19,w19} for the 3-D steady Euler system. The deformation-curl decomposition for the 3-D steady Euler system had played crucial roles in the structural stability of cylindrical transonic shock under three dimensional perturbations of supersonic incoming flows and the exit pressure (See \cite{wx23a}). It was also utilized in \cite{wx23b} to establish the existence and uniqueness of Meyer type smooth transonic flows with nonzero vorticity and positive acceleration to the quasi 2-D Euler flow model in De Laval nozzles. Very recently, the authors in \cite{wz24} proved the existence and uniqueness of smooth axially-symmetric transonic flows to the steady Euler equations with an external force. This had been generalized to the generic three dimensional case in \cite{wx24} and smooth transonic irrotational flows and transonic Beltrami flows with nonconstant proportionality have been constructed in finitely long nozzles with arbitrary cross section.

Motivated by \cite{wxy21b}, we first studies the structural stability of smooth irrotational transonic self-similar flows on the radial distance and the polar angle coordinates. We first analyze the linear mixed type second order equation which is obtained by linearizing the potential equation near the smooth self-similar irrotational transonic flows. By exploring some key properties of the self-similar transonic solutions, we are able to find a multiplier and identify a class of admissible boundary conditions for the linearized mixed type second-order equation. However, as shown in \eqref{m19}, the admissible boundary conditions obtained are oblique, we can not apply the Galerkin's method to construct the approximate solutions. Instead, we resort to the theory developed by Friedrichs and establish the existence of generalized $L^2$ solution to the linearized mixed type second-order equation.

Finally, we further investigate the structural stability of smooth irrotational transonic self-similar flows on the polar and azimuthal angles coordinates. We construct a class of smooth transonic flows with nonzero vorticity depending only on the polar and azimuthal angles which is also close to the self-similar irrotational transonic flows by prescribing suitable boundary conditions on the boundaries. The governing equations are hyperbolic-elliptic coupled and an effective decomposition is introduced by mimicking the deformation-curl decomposition developed in \cite{wex19,w19} for the 3-D steady Euler system. It is worthy to point out that the second elliptic equation satisfied by the radial velocity has a zeroth order term with a positive coefficient, thus the maximum principle and the Hopf's lemma may fail for the radial velocity. Nevertheless, we can show that a uniqueness result for the linearized second order elliptic equation for the radial velocity supplemented with a class of oblique boundary conditions.

The structure of this paper will be organized as follows. In section \ref{self-similar}, we study the self-similar solutions to the Euler system \eqref{SCEQF}. Section \ref{transonic1} will concentrate on the structural stability of the smooth self-similar irrotational transonic flows on the radial distance and polar angle coordinates. In section \ref{transonic2}, we prove the existence and uniqueness of a class of smooth transonic flows with nonzero vorticity which depends only on the polar and azimuthal angles in spherical coordinates.

\section{The self-similar flows to the 3D steady Euler system past an
infinitely long circular cone}\label{self-similar}

\subsection{Smooth self-similar flows to the steady 3D Euler system}

We will analyze \eqref{SCEQF1} with two different types of boundary conditions. Firstly, we will focus on the initial value problem, where we study the behavior of \eqref{SCEQF1} with initial data given at $\theta=\theta_0\in(0,\frac{\pi}2)$.

{\bf Problem I.} Construct a smooth self-similar flows to the steady Euler system near the conic surface $\theta=\theta_0\in (0,\frac{\pi}{2}]$. That is, we are going to solve the following initial value problem to \eqref{SCEQF1} with the initial data at $\theta=\theta_0$
\begin{align}\label{ProblemV}
(\rho,U_1,U_2,U_3,A)(\theta_0)=(\rho_0>0,U_{01}\geq 0,U_{02}\neq 0,U_{03},A_0).
\end{align}
By solving this initial value problem, we will obtain a solution that captures the self-similar behavior of the flow near the conic surface. The solution will depend on the specific choice of initial data and the form of the steady Euler system \eqref{SCEQF1}.

Before stating the main theorem for {\bf Problem I}, it is necessary to introduce some notations. We denote the Mach number in each direction as $M_i(\theta) = \frac{U_i(\theta)}{c(\rho,A)(\theta)}$ for $i = 1, 2, 3$. Let ${\bf M} = (M_1, M_2, M_3)^T$, and define the following three important quantities:
\begin{eqnarray}\label{g1}
g_1(\theta)&=&U_1(\theta)\sin\theta+U_2(\theta)\cos\theta,\\\label{g2}
g_2(\theta)&=&U_2g_1(\theta)+U_3^2\cos\theta=U_1U_2\sin\theta+(U_2^2+U_3^2)\cos\theta,\\\label{g3}
g_3(\theta)&=&U_1(\theta)\cos\theta-U_2(\theta)\sin\theta.
\end{eqnarray}
Note that the ODE system \eqref{SCEQF1} becomes singular when $U_2(\theta) = 0$. Therefore, for our analysis, we will consider only the flow where $U_2(\theta)\neq 0$. In this particular case, the system \eqref{SCEQF1} can be equivalently expressed as follows:
\begin{align}\no
\begin{cases}
\frac{d}{d\theta}(\rho U_2\sin\theta)+2\rho U_1\sin\theta=0,\\
B(\theta)\equiv B_0=\frac12(U_{01}^2+U_{02}^2+U_{03}^2)+\frac{\gamma}{\gamma-1}A_0\rho_0^{\ga-1},\\
U'_1(\theta)=\frac{U_2^2+U_3^2}{U_2},\\
U'_3(\theta)=-(\frac{U_1}{U_2}+\frac1{\tan\theta})U_3,\\
A(\theta)\equiv A_0,\\
(U_1,U_2,U_3)(\theta_0)=(U_{01}\geq 0,U_{02}\neq 0,U_{03}).
\end{cases}
\end{align}
Moreover, if $1-M_2^2(\theta)\neq 0$, we can further conclude that
\begin{align}\label{ODE0}
\begin{cases}
U'_2(\theta)=-\left(U_1+\frac{g_2(\theta)}{(1-M_2^2)U_2\sin\theta}-\frac{U_3^2}{U_2\tan\theta}\right),\\
U'_1(\theta)=\frac{U_2^2+U_3^2}{U_2},\\
U'_3(\theta)=-\frac{g_1U_3}{U_2\sin\theta},\\
B(\theta)\equiv B_0,\ \ \ A(\theta)\equiv A_0,\\
(U_1,U_2,U_3)(\theta_0)=(U_{01}\geq 0,U_{02}\neq 0,U_{03}),
\end{cases}
\end{align}
and
\begin{eqnarray}\label{rhopr}
\rho'(\theta)&=&\left(\frac{M_2^2 g_1(\theta)}{(1-M_2^2)U_2\sin\theta}+\frac{M_3^2}{(1-M_2^2)\tan\theta}\right)\rho\\\no
&=&\frac{g_2(\theta)}{(c^2(\rho,A)-U_2^2)\sin\theta}\rho,\\\label{g1p}
g_1'(\theta)&=&\frac{M_3^2(\sin^2\theta-M_2^2)}{M_2^2(1-M_2^2)\sin\theta}U_2-\frac{1}{(1-M_2^2)\tan\theta}g_1(\theta)\\\no
&=&-\frac{g_2(\theta)\cos\theta}{(1-M_2^2)U_2\sin\theta}+\frac{U_3^2}{U_2\sin\theta},\\\label{g2p}
g_2'(\theta)&=&-\frac{(2-M_2^2)U_1\sin\theta+2U_2\cos\theta}{(1-M_2^2)\sin\theta U_2}g_2(\theta).
\end{eqnarray}
In addition, we have
\be\nonumber
&&(M_1^2)'(\theta)=\frac{2M_1}{M_2}(M_2^2+M_3^2)-(\ga-1)\frac{M_1^2}{1-M_2^2}(\frac{M_2g_1}{c\sin\theta}+\frac{M_3^2}{\tan\theta})\\\label{M1p}
&&=\frac{2M_1}{M_2}(M_2^2+M_3^2)-(\ga-1)\fr{M_1^2\rho'}{\rho},\\\no
&&(M_2^2)'(\theta)=-2M_1M_2-\frac{2M_2 g_1}{(1-M_2^2)c\sin\theta}(1+\frac{\ga-1}2 M_2^2)-\frac{(\ga+1)M_2^2M_3^2}{(1-M_2^2)\tan\theta}\\\label{M2p}
&&=-2M_1M_2-2M_2\frac{g_1}{c\sin\theta}-(\ga+1)M_2^2\frac{\rho'}{\rho},\\\no
&&(M_3^2)'(\theta)=-\frac{2M_3^2g_1}{M_2c\sin\theta}(1+\frac{\ga-1}2\frac{M_2^2}{1-M_2^2})-\frac{(\ga-1)M_3^4}{(1-M_2^2)\tan\theta}\\\label{M3p}
&&=-\frac{2M_3^2g_1}{U_2\sin\theta}-(\ga-1)M_3^2\dfr{\rho'}{\rho},\\\no
&&(|{\bf M}|^2)'(\theta)=-\frac{2M_2g_1}{(1-M_2^2)c\sin\theta}(1+\frac{\ga-1}2 |\M|^2)\\\no
&&\quad-M_3^2\big(\frac2{\tan\theta}+\frac1{(1-M_2^2)\tan\theta}\big((\ga+1)M_2^2+(\ga-1)(M_1^2+M_3^2)\big)\big)\\\label{Mp}
&&=-\frac{2g_2}{(1-M_2^2)c^2\sin\theta}(1+\frac{\ga-1}2 |{\bf M}|^2).
\ee
Therefore, the quantity $g_2(\theta)$ plays a key role in determining the signs of $\rho'(\theta)$ and $\frac{d}{d\theta}|{\bf M}|^2$.

Now, we are ready to state the main theorem for {\bf Problem I}:
\bt\label{ProbITH}
{\it The system \eqref{SCEQF1} with the initial data \eqref{ProblemV} has a unique solution near $\theta=\theta_0$. To describe the behavior of the solution, we divide it into the following cases:
\begin{enumerate}[(1)]
  \item If $1-M_2^2(\theta_0)<0$, then the flow is always supersonic until $1-M_2^2(\theta)$ reduces to $0$ where a singularity may happen.
  \item If $g_2(\theta_0)=0$, then there exist a unique $\theta_{min}\in(0,\theta_0)$ such that Problem I is uniquely solvable in  $\theta\in(\theta_{min},\pi-\theta_{min})$. In this case, $\rho(\theta)\equiv \rho_0, A(\theta)\equiv A_0$ and the velocity has the following formula.
  \begin{enumerate}[(a)]
  \item If $\theta_0\in(0,\frac{\pi}2)$, then $U_{01}>0, U_{02}<0$, $\theta_{min}=\arccos\big(\frac{q_0}{U_{01}}\cos\theta_0\big)\in [0,\theta_0)$ and
  \begin{eqnarray}\label{solu1}
  \begin{cases}
  U_1(\theta)=\frac{q_0^2\cos\theta}{U_{01}\cos\theta_0-U_{02}\sin\theta_0}=\frac{U_{01}}{\cos\theta_0}\cos\theta,\\
  U_2(\theta)=\frac{U_{01}^2\cos^2\theta-q_0^2\cos^2\theta_0}{U_{01}\cos\theta_0\sin\theta},\\
  U_3^2(\theta)=\frac{(q_0^2\cos^2\theta_0-U_{01}^2)(U_{01}^2\cos^2\theta-q_0^2\cos^2\theta_0)}{U_{01}^2\cos^2\theta_0\sin^2\theta},\\
  \end{cases}
  \end{eqnarray}
  where $q_0^2=U_{01}^2+U_{02}^2+U_{03}^3$ and $U_3$ has the same sign as $U_{03}$ and does not change its sign when $\theta$ varies.

  \item If $\theta_0=\frac\pi2$, then $U_{01}=0$, $\theta_{min}=\arcsin\frac{|U_{03}|}{q_0}\in(0,\theta_0)$ and
  \begin{eqnarray}\label{solu2}
  \begin{cases}
  U_1(\theta)=\frac{q_0^2\cos\theta}{-U_{02}},\\
  U_2(\theta)=\frac{U_{02}^2-q_0^2\cos^2\theta}{U_{02}\sin\theta},\\
  U_3^2(\theta)=\frac{(q_0^2-U_{02}^2)(U_{02}^2-q_0^2\cos^2\theta)}{U_{02}^2\sin^2\theta},
  \end{cases}
  \end{eqnarray}
  where $U_3$ has the same sign as $U_{03}$ and does not change its sign when $\theta$ varies.
  \end{enumerate}
  Moreover, $U_2(\th_{min})=U_2(\pi-\th_{min})=U_3(\th_{min})=U_3(\pi-\th_{min})=0$ in both cases $(a)$ and $(b)$.
  \item If $1-M_2^2(\theta_0)>0$, $U_{02}<0$ and $g_2(\theta_0)<0$, then there exist a unique $\theta_*\in (0,\theta_0)$ such that \eqref{ProblemV} is uniquely solvable in $(\theta_*,\theta_0)$ and, for any $\theta\in(\theta_*,\theta_0)$, $\rho'(\theta)>0$, $1-M_2^2(\theta)>0$, and $U_2(\theta_*)=0$, $U_3(\theta_*)=0$.

  \end{enumerate}
}\et

\br{\it
In general, the flows constructed in (2) has non zero vorticity unless $U_{03}=0$ (and thus $U_{01}\sin\theta_0+ U_{02}\cos\theta_0=0$ since $g_2(\th_0)=0$), in which case $U_1(\theta)=\frac{U_{01}}{\cos\theta_0}\cos\theta, U_2(\theta)=-\frac{U_{01}}{\cos\theta_0}\sin\theta$ and $U_3(\theta)\equiv 0$. As a byproduct, since $\rho(\theta)\equiv \rho_0$, the flows in (2) are indeed a class of infinitely smooth incompressible Beltrami flows with a nonconstant proportionality factor. Indeed the velocity field ${\bf u}= U_1(\theta){\bf e}_r+ U_2(\theta){\bf e}_{\theta}+ U_3(\theta){\bf e}_{\varphi}$ presented in (2) satisfies
\be\no
\text{curl }{\bf u}(x)=\frac{U_3(\theta)}{r U_2(\theta)} {\bf u}(x).
\ee
This can be viewed as an interesting supplement to the results proved in \cite{ep16}, which showed that incompressible Beltrami flows with a nonconstant proportionality factor are rare in some sense.
}\er

\br{\it
The flows in (2) can be purely supersonic flows, purely subsonic flows and purely sonic flows with constant Mach number in the whole region.
}\er

\br{\it
The flows in (3) can be purely supersonic flows, purely subsonic flows and smooth transonic flows with decreasing Mach number as $\theta$ increases. The smooth transonic flows constructed here may have a nonzero vorticity, and the sonic points are nonexceptional and noncharacteristic degenerate due to the nontrivial radial and azimuthal components of the velocity $(U_1,U_3)$.
}\er

\bpf[Proof of Theorem \ref{ProbITH}.]
The local unique solvability of \eqref{SCEQF1} with \eqref{ProblemV} and the conclusions in (1) follows directly from the standard theory of ODEs.

To prove the results in (2), assume that $g_2(\theta_0)=0$, then by \eqref{g2p} one has $g_2(\theta)\equiv 0$, which implies that $\rho'(\theta)\equiv 0$ and $\rho(\theta)\equiv \rho_0$.
Since $g_2(\theta)\equiv 0$, the system \eqref{ODE0} reduces to
\begin{eqnarray}\no
\begin{cases}
U'_2(\theta)=-\left(U_1-\frac{U_3^2}{U_2\tan\theta}\right),\\
U'_1(\theta)=U_2+\frac{U_3^2}{U_2},\\
U_3U'_3(\theta)=-\frac{g_1U_3^2}{U_2\sin\theta},\\
B(\theta)\equiv B_0,\ \ \ A(\theta)\equiv A_0,\\
U_1^2+U_2^2+U_3^3\equiv 2(B_0-\frac{\ga A_0\rho_0^{\ga-1}}{\ga-1}).
\end{cases}
\end{eqnarray}
It follows from the first two equations that $(U_2'+U_1)\tan\theta=\frac{U_3^2}{U_2}=U_1'-U_2$, thus $g_3'(\theta)\equiv 0$. Hence ${\bf U}$ should satisfy the following algebraic equations:
\begin{align}\label{s1}
\begin{cases}
U_1U_2\sin\theta+(U_2^2+U_3^3)\cos\theta=0,\\
U_2\sin\theta-U_1\cos\theta\equiv U_{02}\sin\theta_0-U_{01}\cos\theta_0,\\
U_1^2+U_2^2+U_3^2\equiv 2(B_0-\frac{\ga A_0\rho_0^{\ga-1}}{\ga-1})\equiv q_0^2.
\end{cases}
\end{align}

{\bf Case 1.} If $\theta_0\in(0,\frac\pi{2})$, since $g_2(\theta_0)=0$, one has $U_{02}^2+U_{03}^2=-U_{01}U_{02}\tan\theta_0$ and $q_0^2=U_{01}^2-U_{01}U_{02}\tan\theta_0$. Thus $U_{01}>0, U_{02}<0$. Substituting $U_2^2+U_3^2=q_0^2- U_1^2$ into the first equation in \eqref{s1} and using the second equation in \eqref{s1}, one can derive the formula for $U_1$ in \eqref{solu1}. Then the formula for $U_2$ and $U_3$ can be derived directly.

To find the maximal existence interval for $(U_1,U_2,U_3)$ in \eqref{solu1}, we note that
\be\no
U_{01}^2-q_0^2\cos^2\theta_0=U_{01}\sin\theta_0 g_1(\theta_0)=-\frac{U_{01}U_{03}^2}{U_{02}}\sin\theta_0\cos\theta_0\geq 0.
\ee
Thus $U_{01}^2-q_0^2\cos^2\theta_0=0$ if and only if $U_{03}=0$. One may define
\be\no
\theta_{min}=\arccos\bigg(\frac{q_0}{U_{01}}\cos\theta_0\bigg)\in(0,\theta_0),
\ee
then $\theta_{min}\in [0,\theta_0)$ and $\theta_{min}=0$ if and only if $U_{03}=0$. Note that $U_2(\theta_{min})=0$ and $U_2(\pi-\theta_{min})=0$, the maximal existence interval is $(\theta_{min}, \pi-\theta_{min})$.

{\bf Case 2.} If $\theta_0=\frac{\pi}2$, then $U_{01}=0$ and
\begin{align}\nonumber
\begin{cases}
U_1(\theta)=\frac{q_0^2\cos\theta}{-U_{02}},\\
U_2(\theta)=\frac{U_{02}^2-q_0^2\cos^2\theta}{U_{02}\sin\theta},\\
U_3^2(\theta)=\frac{(q_0^2-U_{02}^2)(U_{02}^2-q_0^2\cos^2\theta)}{U_{02}^2\sin^2\theta}.
\end{cases}
\end{align}
The solution exists only if $\theta\in (\theta_{min},\pi-\theta_{min})$, where $\theta_{min}=\arcsin\frac{|U_{03}|}{q_0}$ with $q_0=\sqrt{U_{02}^2+U_{03}^2}$, and $U_2(\theta_{min})=0$ and $U_2(\pi-\theta_{min})=0$.


Since $U_2$ is the velocity in $\theta$-direction, $U_2(\theta_{min})=U_2(\pi-\theta_{min})=0$ means that the flow will not across the surface $\theta=\theta_{min}$ and $\theta=\pi-\theta_{min}$ and will go forward along the surfaces. Moreover, we have $U_3(\theta_{min})=U_3(\pi-\theta_{min})=0$ in both cases. The conclusion in the statement (3) follows directly from the following lemma.
\epf

\bl\label{lem42}
{\it Suppose $g_2(\theta_0)<0$, $1-M_2^2(\theta_0)>0$, $g_1(\theta_0)>0$, $\rho_0>0$, $U_{01}>0$, $U_{02}<0$, $U_{03}^2>0$, there exist a unique $\theta_*\in {\color{red}[}0,\theta_0)$ such that \eqref{ODE0} is uniquely solvable in $(\theta_*,\theta_0)$ and, for any $\theta\in(\theta_*,\theta_0)$,  $g_2(\theta)<0$, $\rho(\theta)>0$, $\rho'(\theta)<0$, $g_1(\theta)>0$, $g_1'(\theta)<0$, $1-M_2^2(\theta)>0$, $(M_2^2)'(\theta)>0$, $(M_3^2)'(\theta)>0$, $U_1(\theta)>0$, $U_2(\theta)<0$, $U_3^2(\theta)>0$, $g_3'(\theta)=\frac{g_2(\theta)}{(1-M_2^2)U_2}>0$ and $U_2(\theta_*)=0$, $U_3(\theta_*)=0$.
}\el

\bpf
Define the set
\begin{align*}
J:=\bigg\{\theta_1\in(0,\theta_0): \eqref{ODE0}\ is\ solvable\ in\ (\theta_1,\theta_0),g_2(\theta)<0 ,  \rho(\theta)>0 , \\ \rho'(\theta)<0 ,  g_1(\theta)>0 , g_1'(\theta)<0,  1-M_2^2(\theta)>0, U_1'(\theta)<0,\\ U_2'(\theta)<0,
(M_2^2)'(\theta)>0 ,  (M_3^2)'(\theta)>0,
U_1(\theta)>0 ,\\  U_2(\theta)<0 , U_3^2(\theta)>0 \ \text{holds for any }\theta\in (\theta_1,\theta_0)\bigg\}.
\end{align*}
By the assumptions on the initial data and the equations \eqref{ODE0}, \eqref{rhopr}, \eqref{g1p}, \eqref{M2p}, \eqref{M3p}, we have
\be\no
g_1'(\theta_0)<0,\ \rho'(\theta_0)<0,\ \frac{d}{d\theta}(M_2^2)(\theta_0)<0,\ \frac{d}{d\theta}(M_3^2)(\theta_0)>0.
\ee
By the classical theory of ODEs, we know that $J$ is a nonempty open set. Define
\begin{align*}
\theta_*=\inf_{\theta_1\in J}\theta_1.
\end{align*}

Case 1: Suppose $\theta_*=0$. Since $U_1'(\theta)<0$ and $U_2'(\theta)<0$ for any $\theta\in (0,\theta_0)$, one can define
\be\no
U_1(0)=\lim_{\theta\to 0}U_1(\theta)\geq 0,\ \ U_2(0)=\lim_{\theta\to 0}U_2(\theta)\leq 0.
\ee
Recall that $g_1(\theta)=U_1(\theta)\sin\theta+U_2(\theta)\cos\theta>0$ and $g_1'(\theta)<0$ for any $\theta\in (0,\theta_0)$, $g_1(0)=U_2(0)=\displaystyle\lim_{\theta\to 0}g_1(\theta)\geq 0$ which yields that $U_2(0)=0$. In addition, $U_3^2(0)=\lim_{\th\to 0}g_2(\th)\leq 0$ implies that $U_3(0)=0$.


Case 2: Suppose $\theta_*>0$, then $U_2(\theta_*)\leq 0$. Assume $U_2(\theta_*)\neq 0$, then $U_2(\theta_*)<0$. By the definition of $\theta_*$, for any small $\epsilon>0$, $1-M_2^2(\theta_*+\epsilon)>0$ and $(M_2^2)'(\theta)>0$ in $(\theta_*,\theta_*+\epsilon)$. Thus, $1-M_2^2(\theta_*)>1-M_2^2(\theta_*+\epsilon)>0$ and thus $\rho(\theta_*)>0$. Since $U_2(\theta)<0$ for any $\theta\in[\theta_*,\theta_0]$, by \eqref{g2p}, we have
\begin{align}\label{g2}
g_2(\theta_*)=g_2(\theta_0)\exp\bigg\{\int_{\theta_0}^{\theta_*}-\frac{(2-M_2^2)U_1\sin\theta+2U_2\cos\theta}{(1-M_2^2)\sin\theta U_2}d\theta\bigg\}<0.
\end{align}
Thus $U_1(\theta_*)>0$ and, by \eqref{rhopr}, \eqref{g1p}, $\rho'(\theta_*)<0$, $g_1'(\theta_*)>0$. Therefore, $g_1(\theta_*)>0$ and by \eqref{M2p}, \eqref{M3p}, $(M_2^2)'(\theta)>0$,  $(M_3^2)'(\theta)>0$. So, by the definition of $J$, we have $\theta_*-\epsilon_0\in J$ for sufficiently small $\epsilon_0>0$, which contradicts to the definition of $\theta_*$. Therefore, $U_2(\theta_*)=0$ and $g_2(\theta_*)=U_3^2(\theta_*)\cos\theta_*\geq 0$. By \eqref{g2}, we have $g_2(\theta_*)\leq 0$. Thus, $g_2(\theta_*)=0$ and $U_3(\theta_*)=0$.
\epf

\br\label{VI43}
{\it By the same argument we can prove the case $U_{03}^2=0$. That is, suppose $g_2(\theta_0)<0$, $1-M_2^2(\theta_0)>0$, $g_1(\theta_0)>0$, $\rho_0>0$, $U_{01}>0$, $U_{02}<0$, $U_{03}^2=0$, there exist a unique $\theta_*\in (0,\theta_0)$ such that \eqref{ODE0} is uniquely solvable in $(\theta_*,\theta_0)$ and, for any $\theta\in(\theta_*,\theta_0)$,  $g_2(\theta)<0$, $\rho(\theta)>0$, $\rho'(\theta)<0$, $g_1(\theta)>0$, $g_1'(\theta)<0$, $1-M_2^2(\theta)>0$, $(M_2^2)'(\theta)>0$, $U_3(\theta)\equiv 0$, $U_1(\theta)>0$, $U_2(\theta)<0$, $g_3'(\theta)=(U_1\cos\theta-U_2\sin\theta)'(\theta)=\frac{g_2(\theta)}{(1-M_2^2)U_2}>0$ and $U_2(\theta_*)=0$.
}\er


By Lemma \ref{lem42} and the above remark, we finish the proof of Theorem \ref{ProbITH}.

\subsection{Supersonic flows passing through an infinitely long conic cone}

In accordance with the findings in \cite{cf48}, when a supersonic constant incoming flow encounters an infinitely long circular cone, a conic shock will develop, connected to the cone's tip, provided that the cone's opening angle is smaller than a critical value. In this context, we will now examine a more intricate scenario: the interaction of a supersonic self-similar incoming flow with a non-zero azimuthal velocity with an infinitely long circular cone. Our objective is to determine the conditions under which a supersonic or transonic shock solution exists, attached to the cone's tip.

Consider the problem of finding a piecewise smooth solution to the system \eqref{SCEQF} on a domain $\Omega_1$ with a discontinuity on a surface $S_f$ parameterized by $(r,\theta=f(r,\varphi),\varphi)$. We define two domains as $\Omega_1^+:=\{(r,\theta,\varphi)\in \Omega_1:\theta>f(r,\varphi)\}$ and $\Omega_1^-:=\{(r,\theta,\varphi)\in \Omega_1:\theta<f(r,\varphi)\}$.
A piecewise smooth solution $(\rho^\pm,U_{1}^\pm,U_{2}^\pm,U_{3}^\pm,A^\pm)\in C^1(\Omega_1^{\pm})$ with a jump on the surface $S_f$ is considered a weak solution to the system \eqref{SCEQF} if it satisfies the system in the respective domains $\Omega_1^{\pm}$ and satisfies the Rankine-Hugoniot jump conditions on the surface $S_f$. Therefore, our goal is to solve:

{\bf Problem II.} To construct a piecewise smooth self-similar solution to \eqref{comeuler3d} in $I=[\theta_*,\frac{\pi}{2}]$,
\be\no
(\rho^\pm(\theta),U_{1}^\pm(\theta),U_{2}^\pm(\theta),U_{3}^\pm(\theta),A^\pm(\theta)),
\ee
which satisfy the following ODE system in $I^-=(\theta_b,\frac{\pi}{2}]$ and $I^+=(\theta_*,\theta_b)$ respectively,
\begin{align}\label{ProbVI}
\begin{cases}
(\rho U_2)'(\theta)+2 \rho U_1+\frac1{\tan\theta}\rho U_2=0,\\
U_2 U_1'-(U^2_2+U_3^2)=0,\\
U_2U_2'+\frac1\rho \frac{d}{d\theta}p(\rho)+U_1U_2-\frac{U_3^2}{\tan\theta}=0,\\
U_2U_3'+U_1U_3+\frac{U_2U_3}{\tan\theta}=0,\\
U_2 A'=0,\\
\end{cases}
\end{align}
with the boundary conditions on both $\theta=\frac{\pi}{2}$ and $\theta=\theta_*$:
\be\label{var0}
&&(\rho,U_1,U_2,U_3,A)(\frac{\pi}{2})=(\rho_0>0,0,U_{02}<0,U_{03},A_0),\\\label{var1}
&&U_2(\theta_*)=0.
\ee
The discontinuity occurs at $\theta=\theta_b\in (\theta_*,\frac{\pi}{2})$ which is unknown and should be determined together with the solutions simultaneously. Across the discontinuity $\theta=\theta_b$, the physical entropy condition $[p]>0$ and the following Rankine-Hugoniot conditions hold,
\begin{align}\label{rh50}
\begin{cases}
[\rho U_2]=0,\\
[\rho U_1U_2]=0,\\
[\rho U_2^2+p]=0,\\
[\rho U_2U_3]=0,\\
[B]=0.
\end{cases}
\end{align}

Define $q_0^2=U_{02}^2+U_{03}^2, c_0=\sqrt{A_0\gamma \rho_0^{\gamma-1}}, c_*^2=\frac{\gamma-1}{\gamma+1}q_0^2+\frac{2}{\gamma+1}c_0^2$ and $a_1^2:=\frac{U_{02}^2}{c_0^2}$.

\bt\label{ProbIITH}
{\it Suppose that $a_1^2>1$ and, for all $t\in(0,1)$,
\be\label{290}
(\ga-1)q_0^4t^2+\bigg(2c_0^2U_{02}^2-(3\ga-1)U_{03}^2q_0^2\bigg)t+2\ga U_{03}^4>0.
\ee
Then, there exists a set $D\subset(0,\frac{\pi}2)$ such that for any $\theta_*\in D$, {\bf Problem II} has at least one conic shock solution $ (\rho^\pm(\theta),U_{1}^\pm(\theta),U_{2}^\pm(\theta),U_{3}^\pm(\theta),A^\pm)$ defined on $I^-=[\theta_b, \frac{\pi}{2}]$ and $I^+=(\theta_*,\theta_b]$ with the shock located at $\theta=\theta_b\in (\theta_*,\frac{\pi}{2})$. Moreover, $U_3^+(\theta_*)=0$ and the state just behind the shock can be classified as follows.
\begin{enumerate}
  \item If $a_{\sharp}^2\geq a_1^2$, the state just behind the shock $\theta=\theta_b$ is always supersonic;
  \item If $a_{\sharp}^2< a_1^2$, there exists a unique $\theta_{\sharp}\in (\underline{\theta},\frac{\pi}{2})$ such that if $\theta_b\in (\underline{\theta},\theta_{\sharp})$, the state just behind the shock is supersonic; if $\theta_b=\theta_{\sharp}$, the state just behind the shock is sonic which means that the shock surface and the sonic surface coincide; if $\theta_b\in (\theta_{\sharp}, \frac{\pi}{2})$, the state just behind the shock is subsonic.
\end{enumerate}
Here $a_{\sharp}^2$ is defined in \eqref{a-sharp}, $\underline{\theta}, \theta_{\sharp}$ are defined in \eqref{under-theta},\eqref{theta-sharp} respectively; all of them depend only on $\rho_0$, $U_{02}$, $U_{03}$, $A_0$ and $\gamma$.

}\et

\br
{\it
Any of the following conditions can imply the validity of \eqref{290}:
\begin{enumerate}[(a)]
  \item $2c_0^2U_{02}^2>(3\ga-1-2\sqrt{2\ga(\ga-1)})q_0^2U_{03}^2$;
  \item $(3\gamma-1)U_{03}^2 q_0^2\geq 2c_0^2U_{02}^2+2(\gamma-1)q_0^4$ and $(\gamma-1)q_0^4+2c_0^2 U_{02}^2+ 2\gamma U_{03}^4>(3\gamma-1)q_0^2 U_{03}^2$.
\end{enumerate}
}\er
\br\label{r21}
{\it We fail to give a precise characterization of the set $D$, therefore for supersonic incoming flows satisfying the conditions in Theorem \ref{ProbIITH} past a given infinitely long circular cone, we do not know when does there exist a conic shock solution and how many conic shock solutions does it have.
}\er

\br\label{r22}
{\it In the case $U_{03}=0$, one can verify that $a_{\sharp}^2<a_1^2$ and we can give a more detailed description about the flow state in downstream. Please refer to the proof of Theorem \ref{ProbIITH}.
}\er



\bpf
First, we consider the specific case where $U_{03}=0$, and it can be observed that the swirl velocity component $U_3$ remains identically zero throughout the flow. Indeed, in the presence of the boundary condition $U_{01} = 0$, as stated in \eqref{var0}, we can conclude that $g_2\left(\frac{\pi}{2}\right) = 0$ based on the definition of $g_2$ given in \eqref{g2}. By Theorem \ref{ProbITH}, we have
\begin{eqnarray}\no
U_1^-(\theta)= q_0\cos\theta,\ \ U_2^-(\theta)= -q_0\sin\theta, \ \ U_3(\theta)\equiv 0.
\end{eqnarray}
where $q_0=|U_{02}|$. Thus we have a horizontal supersonic incoming flow ${\bf u}\equiv (0,0, q_0)^\top$ in Cartesion coordinate and $\rho\equiv \rho_0, A\equiv A_0$ with $A_0\gamma \rho_0^{\gamma-1}>q_0^2$ hits an infinite long symmetric cone, a shock must occur due to the sharp change of direction of the velocity and the compressibility of the flow. This phenomena had been investigated ninety years ago by \cite{cf48,tm33}. We should state clearly that most of the results obtained in the following for this special case are not new. The cylindrical coordinates were used in \cite{cxy02,cf48,xy08} to solve this problem, the shock was treated as an oblique shock problem and the shock polar and the apple curve were introduced to illustrate when does there exists a conic shock solution. Different from previous results, we use the spherical coordinates to describe this phenomena where the shock can be treated as a normal shock.

It follows from the Rankine-Hugoniot condition \eqref{rh50} that $U_1$ and $B$ have no jump across the shock, and $U_2$ will jump from a relative supersonic state to a relative subsonic state. Suppose there is a shock occurring at $\theta=\theta_b$, then
\begin{eqnarray}\nonumber
U_1^+(\theta_b)=U_1^-(\theta_b)=q_0\cos\theta_b,\ \ \ U_3^+(\theta_b)=0,\ \ \ B^+(\theta_b)=B_0.
\end{eqnarray}
It follows from the first and third equations in \eqref{rh50} that $[U_2+\frac{c^2}{\gamma U_2}]=0$, which implies that $[U_2+ \frac{K(\theta)}{U_2}]=0$, where
\be\nonumber
K(\theta):=\frac{2}{\ga+1}c^2(\theta)+\frac{\ga-1}{\ga+1}U_2^2(\theta)\equiv\frac{(\ga-1)}{\ga+1}(2B_0-U_1^2(\theta)).
\ee
Note that $K(\theta)$ has no jump when crossing the shock $\theta=\theta_b$: $[K(\theta_b)]=0$. Thus the Prandtl's relation holds
\be\label{prandtl}
U_2^+(\theta_b) U_2^-(\theta_b)= K(\theta_b).
\ee
To guarantee that $|U_2^-(\theta_b)|> |U_2^+(\theta_b)|$, by \eqref{prandtl}, it is equivalent to
\be\no
|U_2^-(\theta_b)|^2> K(\theta_b)=\frac{2}{\ga+1}(c^-(\theta))^2+\frac{\ga-1}{\ga+1}|U_2^-(\theta)|^2.
\ee
Thus it is required that $q_0^2\sin^2\theta_b> A_0\gamma \rho_0^{\gamma-1}(:=c_0^2)$, i.e.
\begin{eqnarray}\nonumber
\underline{\theta}:=\arcsin\frac{1}{M_0}<\theta_b<\frac{\pi}{2}, \ \ M_0:=\frac{q_0}{c_0}>1.
\end{eqnarray}
From \eqref{prandtl} we deduce that
\begin{align}\no
\begin{cases}
U_2^+(\theta_b)=\frac{K(\theta_b)}{U_2^-(\theta_b)}=\frac{K(\theta_b)}{-q_0\sin\theta_b}<0,\\
\rho^+=\frac{\rho^-U_2^-}{U_2^+}=\frac{\rho_0q_{0}^2\sin^2\theta_b}{K(\theta_b)},\\
A^+_0=(\frac{K(\theta_b)}{\rho_0 q_{0}^2\sin^2\theta_b})^{\gamma-1}\frac{\gamma-1}{\gamma}(B_0-\frac12 q_{0}^2\cos^2\theta_b-\frac12\frac{K^2(\theta_b)}{q_0^2\sin^2\theta_b}).
\end{cases}
\end{align}
Next, we introduce the shock polar in the $(U_1^+, U_2^+)$ plane. Since
\be\no
K(\theta_b)=\frac{\gamma-1}{\gamma+1}(q_0^2-(U_1^+(\theta_b))^2)+ \frac{2c_0^2}{\gamma+1},
\ee
one can represent $U_2^+(\theta_b)$ as a function of $U_1^+(\theta_b)$: $U_2^+(\theta_b)=G(U_1^+(\theta_b))$ (see Figure 1), where
\begin{eqnarray}\no
G(s)=-\frac{1}{\gamma+1} \frac{1}{\sqrt{q_0^2-s^2}}\big((\gamma-1)(q_0^2-s^2)+2 c_0^2\big).
\end{eqnarray}
Note that, for $0<s< q_0 \cos \underline{\theta}$,
\be\no
G'(s)=\frac{1}{\gamma+1} s (q_0^2-s^2)^{-\frac{3}{2}} (q_0^2-s^2+ 2c_0^2)>0.
\ee
Furthermore, we can deduce the following information:
\be\no
g_2^+(\theta_b)&=&U_1^+(\theta_b)U_2^+(\theta_b)\sin\theta_b+(U_2^+(\theta_b))^2\cos\theta_b\\\no
&=&U_2^+(\theta_b)\cos\theta_b(U_1^+(\theta_b)\tan\theta_b+U_2^+(\theta_b))\\\no
&=&U_2^+\cos\theta_b(U_2^+-U_2^-)(\theta_b)<0,
\ee
and $1-(M_2^+)^2(\theta_b)>0$, $g_1^+(\theta_b)=\cos\theta_b(U_2^+-U_2^-)(\theta_b)>0$, $\rho_0^+(\theta_b)>0$, $U_{1}^+(\theta_b)>0$, $U_{2}^+(\theta_b)>0$, $U_{3}^+(\theta_b)=0$. Therefore, by Remark \ref{VI43}, there exists a unique $\theta_*\in(0,\theta_b)$  depending on the the incoming flow and the position of the shock $\theta_b$ such that there exists a unique smooth solution
$(\rho^+(\theta),U_{1}^+(\theta),U_{2}^+(\theta),0,A^+(\theta))$ to \eqref{ProbVI} in $(\theta_*,\theta_b]$ with initial data
$(\rho^+(\theta_b),U_{1}^+(\theta_b),U_{2}^+(\theta_b),0,A^+(\theta_b))$. Moreover, we have $U_2^+(\theta_*)=0$. In this case, the constant $\theta_*$ can be determined as the zero point of $(U_1^+)'(\theta_*)=0$, where $U_1^+(\theta)$ satisfies a second-order differential equation
\be\no
&&U_1''(\theta)+ U_1(\theta)+\frac{[(\gamma-1)B_0-\frac{\gamma-1}{2}(U_1^2+(U_1')^2)](U_1\sin\theta-U_1'\cos\theta)}{[(\gamma-1)B_0-\frac{\gamma-1}{2} U_1^2-\frac{\gamma+1}{2}(U_1')^2]\sin\theta}=0,\\\no
&&U_1(\theta_b)=U_1^+(\theta_b),\ \ U_1'(\theta_b)=U_2^+(\theta_b).
\ee

Now we clarify the state behind the shock $\theta=\theta_b$:
\be\no
&\quad&(U_1^+)^2(\theta_b)+(U_2^+)^2(\theta_b)-c^2(\rho^+(\theta_b),A_0^+)\\\no
&=&\frac{\gamma+1}{2}(q_0^2 \cos^2\theta_b+ \frac{K^2(\theta_b)}{q_0^2\sin^2\theta_b})-(\frac{\gamma-1}2 q_0^2+ c_0^2)\\\no
&=&\frac{q_0^2}{t M_0^4} h_0(M_0^2 t),
\ee
where $t=\sin^2\theta_b\in(0,1)$ and
\be\label{h1}
h_0(\tau):=-\frac{2\ga}{\ga+1}\tau^2+(\frac{\ga-3}{\ga+1}+M_0^2)\tau+\frac2{\ga+1}.
\ee
Note that for the function $h_0(\tau)$, we have $h_0(1)=M_0^2-1>0$ and $h_0(M_0^2)=-\frac{1}{\gamma+1} ((\gamma-1)M_0^2+2)(M_0^2-1)<0$.
Thus, it can be verified that the algebraic equation $h_0(\tau) = 0$ has a unique positive root denoted by $\tau = a_{\sharp}^2>1$, where
\be\label{a-sharp}
a_{\sharp}^2=\frac{1}{4\ga}((\ga-3)+(\ga+1)M_0^2+\sqrt{(\ga-3+(\ga+1)M_0^2)^2+16\ga})<M_0^2.
\ee
Therefore we can define $\theta_{\sharp}=\arcsin \frac{|a_{\sharp}|}{M_0}\in (\underline{\theta},\frac{\pi}{2})$. If $\theta_b\in (\underline{\theta},\theta_{\sharp})$, the state just behind the shock is supersonic; if $\theta_b\in (\theta_{\sharp},\frac{\pi}{2})$, the state just behind the shock is subsonic; if $\theta=\theta_{\sharp}$, the state just behind the shock is sonic, which means that the shock front and the sonic surface coincide. By the standard theory of ODEs, the integral curves (for instance, the curves $L_0$, $L_1$ and $L_2$ in Figure 1) starting from any point $P$ at the shock polar on the $(U_1^+,U_2^+)$ plane will tend to a point $Q(q,0)$ between $(\widetilde{q}_0,0)$ and $(q_0,0)$, where
\be\no
\widetilde{q}_0=-G(0)=\frac{1}{(\gamma+1)q_0}((\gamma-1)q_0^2+ 2c_0^2).
\ee
Therefore the radial velocity on the cone surface has a lower and upper bound and is monotonic decreasing when the shock position $\theta=\theta_b$ increases. This monotonicity relation seems to be new and has not been found before.

\begin{center}
\includegraphics[height=8cm, width=10cm]{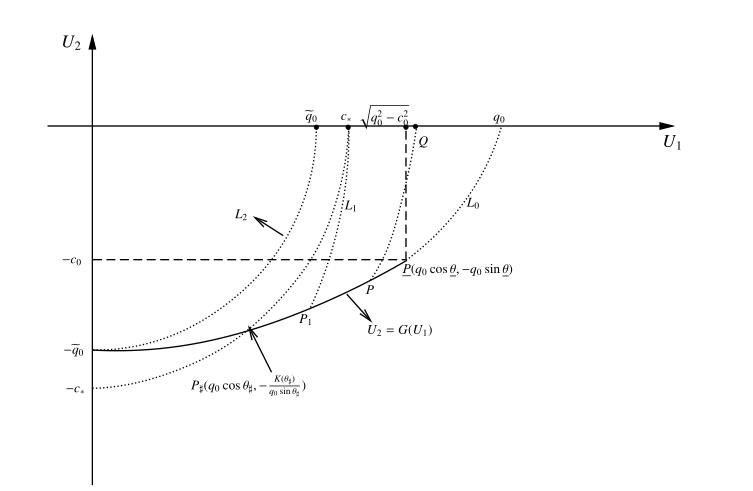}\\
{\small Figure 1: shock polar in $U_1-U_2$ plane.}
\end{center}
Since $(U_1^+)'(\theta)<0, (U_2^+)'(\theta)<0$ and
\be\no
(|{\bf M}|^2)'(\theta)=-\frac{2g_2}{(1-M_2^2)c^2\sin\theta}(1+\fr{\ga-1}2 |{\bf M}|^2)>0,
\ee
for any $\theta\in (\theta_*,\theta_b)$, the integral curve starting from $P_{\sharp}=(q_0\cos\theta_{\sharp},\frac{K(\theta_{\sharp})}{-q_0\sin\theta_{\sharp}})$ will move into the circle $\{(U_1,U_2):U_1^2+U_2^2=c_*^2\}$. There also exists a unique $\theta_s\in (\underline{\theta},\theta_{\sharp})$ such that the integral curve starting from $P_1=(q_0\cos\theta_s,\frac{K(\theta_s)}{-q_0\sin\theta_s})$ will tend to $(c_*,0)$ on the $(U_1,U_2)$ plane. We can conclude that
\begin{enumerate}[(1)]
  \item if $\theta_b\in (\underline{\theta},\theta_s)$, there exists a supersonic-supersonic shock and the flow is uniformly supersonic in $I^+$;

  \item if $\theta_b=\theta_s$, there exists a supersonic-supersonic shock and the flow is supersonic in $I^+$ but degenerates to sonic on the cone surface;

  \item if $\theta_b\in (\theta_s,\theta_{\sharp})$, there exists a supersonic-supersonic shock and the flow changes smoothly from supersonic to subsonic in $I^+$;

  \item if $\theta_b\in (\theta_{\sharp},\frac{\pi}{2})$, there exists a supersonic-subsonic shock and the flow is uniformly subsonic in $I^+$;

  \item if $\theta=\theta_{\sharp}$, there exists a supersonic-sonic shock and the flow is subsonic in $I^+$ but degenerates to sonic at the shock front.

\end{enumerate}


Let us now consider the general case where the swirl velocity $U_{03}$ is non-zero. In this case, we have a supersonic incoming flow with $M_{02}^2>1$ and $U_{03}^2>0$. By the boundary condition $U_{01} = 0$, as stated in \eqref{var0}, we can conclude again that $g_2\left(\frac{\pi}{2}\right) = 0$ based on the definition of $g_2$ given in \eqref{g2}. Consequently, by Theorem \ref{ProbITH}, we have
\begin{align}\no
\begin{cases}
U_1^-(\theta)=\dfr{q_0^2\cos\theta}{-U_{02}},\\
U_2^-(\theta)=\dfr{U_{02}^2-q_0^2\cos^2\theta}{U_{02}\sin\theta},\\
(U_3^-)^2(\theta)=q_0^2-((U_{1}^-)^2+(U_{2}^-)^2),\\
\rho^-(\theta)\equiv\rho_0,\ \ \ A^-(\theta)\equiv A_0,
\end{cases}
\end{align}
for $\theta\in (\theta_{min},\frac{\pi}{2}]$ , where $\theta_{min}=\arcsin\frac{|U_{03}|}{q_0}$ and $q_0^2=U_{02}^2+U_{03}^2$.

Suppose there is a shock occurs at $\theta=\theta_b>\theta_{min}$, by the Rankine-Hugoniot condition \eqref{rh50}, we have
\begin{eqnarray}\label{VI2}
\begin{cases}
U_1^+(\theta_b)=\frac{q_0^2\cos\theta_b}{-U_{02}}>0,\ \ U_2^+(\theta_b)=\frac{K(\theta_b)}{U_2^-(\theta_b)}<0,\\
U_3^+(\theta_b)=U_3^-(\theta_b),\ \ \ B=B_0,\\
\rho^+(\theta_b)=\frac{\rho^-U_2^-(\theta_b)}{U_2^+(\theta_b)},\\
A^+_0=(\rho^+(\theta_b))^{-\ga+1}\frac{\ga-1}{\ga}\bigg(B_0-\frac12\bigg( (U_{1}^+)^2+(U_{2}^+)^2+(U_{3}^+)^2\bigg)\bigg).
\end{cases}
\end{eqnarray}
where
\be\no
K(\theta):=\frac{2}{\ga+1}c^2(\theta)+\frac{\ga-1}{\ga+1}U_2^2(\theta)\equiv\frac{2(\ga-1)}{\ga+1}B_0-\frac{\ga-1}{\ga+1}(U_1^2+U_3^2)(\theta).
\ee
The requirement that $|U_2^-(\theta_b)|> |U_2^+(\theta_b)|$ is equivalent to $|U_2^-(\theta_b)|^2>c_0^2$, i.e. $k(\sin^2\theta_b)>0$, where
\begin{eqnarray}\no
k(t)= q_0^4 t^2- (2q_0^2 U_{03}^2 +c_0^2 U_{02}^2) t + U_{03}^4.
\end{eqnarray}
Since $k(\sin^2\theta_{min})=-\frac{c_0^2 U_{02}^2 U_{03}^2}{q_0^2}<0$ and $k(1)=(U_{02}^2-c_0^2)U_{02}^2>0$, there exists a unique $\underline{\theta}\in (\theta_{min},\frac{\pi}{2})$ such that $k(\sin^2\theta_b)>0$ for any $\theta_b\in (\underline{\theta},\frac{\pi}{2})$. Indeed, $\underline{\theta}$ satisfies
\be\label{under-theta}
\frac{U_{03}^2}{q_0^2}<\sin^2\underline{\theta}= \frac{2 q_0^2 U_{03}^2+ c_0^2 U_{02}^2+ c_0|U_{02}|\sqrt{4q_0^2 U_{03}^2+c_0^2 U_{02}^2}}{2q_0^4}<1.
\ee
Next, we evaluate
\begin{align}\no
g_2^+(\theta_b)=&\ U_1^+U_2^+\sin\theta_b+((U_2^+)^2+(U_3^+)^2)\cos\theta_b\\\no
=&\bigg(\dfr{q_0^2\sin\theta_b}{-U_{02}}U_2^++(U_2^+)^2+q_0^2-(U_1^-)^2-(U_2^-)^2\bigg)\cos\theta_b\\\no
=&\bigg(\dfr{q_0^2\sin\theta_b}{U_{02}}(U_2^--U_2^+)+(U_2^+)^2-(U_2^-)^2\bigg)\cos\theta_b\\\no
=&(U_2^+-U_2^-)\cos\theta_b(U_2^+-\dfr{q_0^2-U_{02}^2}{U_{02}\sin\theta_b})\\\no
=&(U_2^+-U_2^-)\cos\theta_b(\dfr{2 c_0^2}{(\ga+1)U_2^-}+\dfr{\ga-1}{\ga+1}U_2^--\dfr{q_0^2-U_{02}^2}{U_{02}\sin\theta_b}),
\end{align}
which is negative if and only if $t:=\sin^2\theta_b\in(0,1)$ satisfies \eqref{290}.

In addition, $1-(M_2^+)^2(\theta_b)>0$, $g_1^+(\theta_b)>g_1^-(\theta_b)=-\fr{U_{03}^2}{U_{02}\sin\theta}\cos\theta>0$, $\rho_0^+(\theta_b)>0$, $U_{1}^+(\theta_b)>0$, $U_{2}^+(\theta_b)<0$, $U_{3}^+(\theta_b)=U_{3}^-(\theta_b)$. Therefore, by Lemma \ref{lem42}, there exists a unique $\theta_*\in [0,\theta_b)$ depending on the the incoming flow and the position of the shock $\theta_b$ such that there exists a unique smooth solution
$(\rho^+(\theta),U_{1}^+(\theta),U_{2}^+(\theta),U_{3}^+(\theta),A^+(\theta))$ to \eqref{ProbVI} in $(\theta_*,\theta_b]$ with initial data $(\rho^+(\theta_b),U_{1}^+(\theta_b),U_{2}^+(\theta_b),U_{3}^+(\theta_b),A^+(\theta_b))$. Moreover, we have $U_2^+(\theta_*)=U_3^+(\theta_*)=0$.


Next, we are going to classify the flow pattern after the shock at $\theta=\theta_b$. By \eqref{VI2}, we have, at $\theta=\theta_b$,
\begin{align*}
&(U_1^+)^2+(U_2^+)^2+(U_3^+)^2-c^2(\rho^+,A_0^+)\\
&=(U_1^+)^2+(U_2^+)^2+(U_3^+)^2-(K_0-\dfr{\ga-1}{\ga+1}(U_2^+)^2)\dfr{\ga+1}2\\
&=(U_1^-)^2+(U_3^-)^2-\dfr{\ga+1}2(K_0-(U_2^+)^2)\\
&=(U_1^-)^2+(U_3^-)^2-\dfr{\ga+1}2 \dfr{K_0}{(U_2^-)^2}((U_2^-)^2-K_0)\\
&=(U_1^-)^2+(U_3^-)^2+\dfr{K_0}{(U_2^-)^2}((c^-)^2-(U_2^-)^2)\\
&=\frac{K(\theta_b)}{(M_2^-)^2}\bigg(1-(M_2^-)^2+\dfr{(M_2^-)^2(M_1^2+M_3^2)}{\dfr2{\ga+1}+\dfr{\ga-1}{\ga+1}(M_2^-)^2}\bigg)\\
&=\frac{K(\theta_b)}{(M_2^-)^2(\dfr2{\ga+1}+\dfr{\ga-1}{\ga+1}(M_2^-)^2)}h_0((M_2^-)^2(t)),
\end{align*}
where $h_0$ is defined as in \eqref{h1}. The function $h_0(\tau)>0$ for $1<\tau<a_{\sharp}^2$ and $h_0(\tau)<0$ for $\tau>a_{\sharp}^2$.

Denote $a_0=\frac{q_0^2}{U_{02}c_0}$, and $t=\sin^2\theta_b$, then $|a_0|>|a_1|>1$ (only when $U_{03}^2=0$, $q_0=-U_{02}$ and $|a_0|=|a_1|$). Then $(M_2^-)^2(\theta_b)=\dfr{(a_1-a_0(1-t))^2}{t}$. 
Since
\begin{align*}
\dfr{d}{dt}(M_2^-)^2(t)=-\dfr{(a_0-a_1)^2}{t^2}+a_0^2,
\end{align*}
then $\dfr{d}{dt}(M_2^-)^2(t)>0$ if $t>1-\frac{U_{02}^2}{q_0^2}=\sin^2 \theta_*$ and $\frac{d}{dt}(M_2^-)^2(t)<0$ if $t<\sin^2 \theta_*$. Since $(M_2^-)^2(\sin^2\underline{\theta})=1$ and $(M_2^-)^2(1)=a_1^2>1$, then
\be\no
1<(M_2^-)^2(\sin^2\theta_b)<a_1^2,\quad \forall \theta_b\in (\underline{\theta}, \frac{\pi}{2}).
\ee

\begin{enumerate}[(i)]
  \item If $a_{\sharp}^2\geq a_1^2$, $(M_2^-)^2(\sin^2\theta_b)<a_{\sharp}^2$ for any $\theta_b\in (\underline{\theta}, \frac{\pi}{2})$, thus $h_0((M_2^-)^2(\sin^2\theta_b))>0$ for any $\theta_b\in (\underline{\theta}, \frac{\pi}{2}))$, which means that the state behind the shock is always supersonic and the flow may smoothly change from supersonic to subsonic in downstream;
  \item If $a_{\sharp}^2< a_1^2$, there exists a unique $\theta_{\sharp}\in (\underline{\theta},\frac{\pi}{2})$ such that \be\label{theta-sharp}
      (M_2^-)^2(\sin^2\theta_{\sharp})=a_{\sharp}^2.
      \ee
      By the monotonicity of $(M_2^-)^2(\sin^2\theta)$ with respect to $\theta$, we find that if $\theta_b\in (\underline{\theta},\theta_{\sharp})$, the state just behind the shock is supersonic and the flow may smoothly change from supersonic to subsonic in downstream; if $\theta_b=\theta_{\sharp}$, the state just behind the shock is sonic which means that the shock surface and the sonic surface coincide; if $\theta_b\in (\theta_{\sharp}, \frac{\pi}{2})$, the state just behind the shock is subsonic and the flow is uniformly subsonic in downstream.
\end{enumerate}
\epf

\section{The stability analysis of smooth self-similar irrotational transonic flows on the radial distance and polar angle plane}\label{transonic1}

Motivated by our prior research on smooth irrotational transonic flows within a concentric cylinder \cite{wxy21a,wxy21b}, we delve into the structural stability analysis of smooth self-similar irrotational transonic flows in this section. Initially, we aim to construct a class of such flows that possess specific properties. Here, we consider the potential flow case, i.e. $U_3\equiv 0$, and the problem \eqref{ODE0} can be simplified to the following set of equations:
\begin{align}\label{s1}
\begin{cases}
\overline{U}_{1}'(\theta)=\overline{U}_{2}(\theta),\\
\overline{U}_{2}'(\theta)=-\overline{U}_{1}(\theta)-\frac{\overline{U}_{1}(\theta)\sin\theta+\overline{U}_{2}(\theta)\cos\theta}{(1-\overline{M}_{2}^2)\sin\theta},\\
c^2(\overline{\rho})=A_0\gamma \overline{\rho}^{\ga-1}=(\ga-1)(B_0-\frac12(\overline{U}_{1}^2+\overline{U}_{2}^2)).
\end{cases}
\end{align}
We supplement equations \eqref{s1} with the initial data given by
\begin{equation}\label{s01}
(\overline{U}_{1},\overline{U}_{2})(\theta_{so})=(U_{01},U_{02}),
\end{equation}
where $\theta_{so}\in (0,\frac{\pi}{2})$, $U_{01}>0$, $U_{02}>0$, and $U_{01}^2+U_{02}^2=\frac{2(\gamma-1)}{\gamma+1}B_0$ in order for the flow $(\overline{U}_1,\overline{U}_2)$ to become sonic at $\theta=\theta_{so}$. It is worth noting that $g_1(\theta)= \overline{U}_{1}(\theta)\sin\theta+\overline{U}_{2}(\theta)\cos\theta$. Consequently, equations \eqref{g1p}, \eqref{M1p}, \eqref{M2p}, and \eqref{Mp} can be simplified as follows:
\begin{eqnarray}\label{s22}
&&g_1'(\theta)=-\frac{g_1(\theta)}{(1-\overline{M}_{2}^2)\tan\theta},\\\label{s3}
&&\frac{d}{d\theta} (|\ol{{\bf M}}|^2)= -\frac{ \overline{M}_{2} g_1(\theta)}{(1-\overline{M}_{2}^2)c(\overline{\rho})\sin\theta} \left(2+(\gamma-1)|\ol{{\bf M}}|^2\right),\\\label{s4}
&&\frac{d}{d\theta} (\overline{M}_{1}^2)=2\overline{M}_{1}\overline{M}_{2} -\frac{(\gamma-1)\overline{M}_{1}^2 \overline{M}_{2} g_1(\theta)}{(1-\overline{M}_{2}^2)c(\overline{\rho})\sin\theta},\\\label{s5}
&&\frac{d}{d\theta} (\overline{M}_{2}^2)=-2\overline{M}_{1}\overline{M}_{2} -\frac{\overline{M}_{2} g_1(\theta)}{(1-\overline{M}_{2}^2)c(\overline{\rho})\sin\theta}\left(2+(\gamma-1) \overline{M}_{2}^2\right).
\end{eqnarray}
In this scenario, we have $0<\overline{M}_{1}(\theta_{so}), \overline{M}_{2}(\theta_{so})<1$, $|\overline{{\bf M}}(\theta_{so})|^2=1$, and $g_1(\theta_{so})>0$. Consequently, based on equation \eqref{s22}, we can infer that $g_1(\theta)>0$ for any $\theta$ in the vicinity of $\theta_{so}$. By using \eqref{s1}, \eqref{s3} and \eqref{s5}, this further implies that near $\theta_{so}$, we have the following:
\begin{eqnarray}\no\begin{cases}
\overline{U}_{1}'(\theta)>0,\ \  \overline{U}_{2}'(\theta)<0,\ \\
\frac{d}{d\theta} (|\ol{{\bf M}}|^2)<0,\ \  \frac{d}{d\theta} (\overline{M}_{2}^2)<0.
\end{cases}\end{eqnarray}
Therefore, based on the previous deductions, we can derive the following proposition regarding the behavior of the background solution near the sonic state:

\begin{proposition}\label{transonic-bg}
{\it There exists an interval $[\theta_-,\theta_+]$ such that $0<\theta_-<\theta_{so}<\theta_+<\frac{\pi}{2}$ and the ODE system \eqref{s1} with the initial data \eqref{s01} has a unique smooth solution $(\overline{U}_{1},\overline{U}_{2})$ on this interval. Moreover, the total Mach number $|\ol{{\bf M}}(\theta)|^2$ is strictly decreasing on $[\theta_-,\theta_+]$ with $|\overline{{\bf M}}(\theta_{so})|^2=1$, and the following properties hold
\be\no
\overline{U}_{1}(\theta)>0,\ \overline{U}_{2}(\theta)>0,\ 1-\overline{M}_{1}^2(\theta)> 0, \ 1-\overline{M}_{2}^2(\theta)> 0, \forall \theta\in [\theta_-,\theta_+].
\ee
}\end{proposition}

We are not concerned with the maximal existence interval of the solution $(\overline{U}_{1}, \overline{U}_{2})$ in Proposition \ref{transonic-bg} at this stage. Our focus now shifts to investigating the structural stability of this smooth irrotational transonic flow within the domain $\Omega=(r_0,r_1)\times (\theta_-,\theta_+)$. We will study the linear mixed-type second-order equations obtained by linearizing the smooth self-similar transonic flows presented in Proposition \ref{transonic-bg}, and identify a class of admissible boundary conditions. 

Note that the corresponding potential function to the self-similar solution $(\overline{U}_{1}, \overline{U}_{2})$ in Proposition \ref{transonic-bg} is $\overline{\Phi}(r,\theta)= r \overline{U}_{1}(\theta)$, so
\begin{eqnarray}\nonumber
\partial_r \overline{{\Phi}}= \overline{U}_{1}(\theta),\ \frac{1}{r}\partial_{\theta}\overline{{\Phi}}= \overline{U}_{1}'(\theta)= \overline{U}_{2}(\theta),\ \partial_{r\theta}^2 \overline{{\Phi}}= \overline{U}_{1}'(\theta).
\end{eqnarray}
Now we consider the irrotational spherical symmetric transonic flow without swirl:
\begin{eqnarray}\label{u_in_r}
{\bf u}= U_1(r,\theta) {\bf e}_r + U_2(r,\theta){\bf e}_{\theta}=\nabla \Phi= \partial_r \Phi {\bf e}_r + \frac{1}{r} \partial_{\theta}\Phi {\bf e}_{\theta},
\end{eqnarray}
where $\Phi=\Phi(r,\theta)$ is a potential function.
It follows from the Bernoulli's law that
\begin{eqnarray}\label{density-formula}
  \rho= H(|\nabla \Phi|^2)= (\gamma-1)^{\frac{1}{\gamma-1}} \big(1-\frac12 |\nabla \Phi|^2\big)^{\frac{1}{\gamma-1}}.
\end{eqnarray}
By substituting expressions \eqref{u_in_r} and \eqref{density-formula} into the density equation (the first equation in \eqref{SCEQF}), we obtain:
\begin{equation*}
  \partial_r(H(|\nabla \Phi|^2)\partial_r\Phi)+\frac1r \partial_{\theta}(H(|\nabla \Phi|^2)\frac1r \partial_{\theta}\Phi)+ \frac2 r H(|\nabla \Phi|^2) \partial_r \Phi + \frac{H(|\nabla \Phi|^2)}{r^2\tan\theta} \partial_{\theta}\Phi=0,
\end{equation*}
which can be rewritten as
\begin{eqnarray}\no
  &&\big[c^2(H)-(\partial_r\Phi)^2\big] \partial_r^2\Phi-\frac{2}{r^2}\partial_r \Phi \partial_{\theta}\Phi \partial_{r\theta}^2\Phi+ \big[c^2(H)-(\frac1r\partial_{\theta}\Phi)^2\big] \frac{1}{r^2}\partial_{\theta}^2\Phi \\\label{eq_3_10}
  &&\quad\quad+ \frac2r c^2(H) \partial_r \Phi + \frac{c^2(H)}{r^2\tan\theta} \partial_{\theta}\Phi + \frac{1}{r^3} \partial_r \Phi (\partial_{\theta}\Phi)^2=0.
\end{eqnarray}
By denoting the difference between $\Phi$ and the background solution $\overline{\Phi}$ as $\Psi=\Phi-\overline{\Phi}$, equation \eqref{eq_3_10} can be rewritten as follows:
\begin{eqnarray}\nonumber
  &&(c^2(H)-(\partial_r\Phi)^2)\partial_r^2 \Psi - \frac{2}{r^2}\partial_r \Phi \partial_{\theta}\Phi \partial_{r\theta}^2\Psi+ \big(c^2(H)-(\frac1r\partial_{\theta}\Phi)^2\big)\frac{1}{r^2}\partial_{\theta}^2\Psi\\\nonumber
  &&+\big[c^2(H)-(\frac1r\partial_{\theta}\Phi)^2-(c^2(\overline{\rho})-(\frac1r\partial_{\theta}\overline{\Phi})^2)\big]\frac{1}{r^2}\partial_{\theta}^2\overline{\Phi}
  +\frac{2}{r}c^2(H) \partial_r \Psi\\\nonumber
  && + \frac{2}{r}(c^2(H)-c^2(\overline{\rho})) \partial_r \overline{\Phi} +\frac{c^2(H)}{r^2\tan\theta} \partial_{\theta} \Psi+ \frac{1}{r^2\tan \theta}(c^2(H)-c^2(\overline{\rho}))\partial_{\theta}\overline{{\Phi}}\\\nonumber
  &&+\frac{1}{r^2}(\partial_{\theta} \Phi)^2 \partial_r \Psi+ \frac{1}{r^2}\big[(\partial_{\theta} \Phi)^2-(\partial_{\theta} \overline{\Phi})^2\big] \partial_r \overline{\Phi}=0,
\end{eqnarray}
which can be expressed as
\begin{equation}\label{s100}
\mathcal{L}\Psi = F(\nabla \Psi)
\end{equation}
where $\mathcal{L}\Psi$ and $F(\nabla \Psi)$ are denoted as follows:
\begin{eqnarray}\label{s100}
\mathcal{L}\Psi&:=&A_{11}(\nabla \Psi) \partial_r^2 \Psi + \frac{2 A_{12}(\nabla \Psi)}{r} \partial_{r\theta}^2 \Psi +\frac{A_{22}(\nabla \Psi)}{r^2} \partial_{\theta}^2 \Psi\\\no
&\quad& + \frac{e_1(\theta)}{r}\partial_r \Psi + \frac{e_2(\theta)}{r^2} \partial_{\theta} \Psi,
\end{eqnarray}
where
\begin{eqnarray}\label{s101}\begin{cases}
A_{11}(\nabla \Psi):= c^2(H)- (\overline{U}_{1}+\partial_r \Psi)^2,\\
A_{12}(\nabla \Psi)=- (\overline{U}_{1}+\partial_r \Psi)(\overline{U}_{2}+\frac{\partial_{\theta} \Psi}r),\\
A_{22}(\nabla \Psi):= c^2(H)- (\overline{U}_{2}+\frac{\partial_{\theta} \Psi}r)^2,\\
e_1(\theta)
:=2 c^2(\overline{\rho}) -\overline{U}_{2}^2+\frac{(\ga-1)\overline{U}_{1}\overline{M}_{2}^2(\overline{U}_{1}\sin\theta+\overline{U}_{2}\cos\theta)}{(1-\overline{M}_{2}^2)\sin\theta}>0,\\
e_2(\theta):=(3-\gamma)\overline{U}_{2}(\overline{U}_{1}+\overline{U}_{2}\cot\th)+\frac{(\ga+1)\overline{U}_{2}}{1-\overline{M}_{2}^2}(\overline{U}_{1}+\overline{U}_{2}\cot\theta)\\
\quad\quad \quad+(c^2(\overline{\rho})-2\overline{U}_{2}^2)\cot\theta,
\end{cases}\end{eqnarray}
and
\begin{eqnarray}\nonumber
&&F(\nabla\Psi):=\frac{2}{r^2}\overline{U}_{1}' \partial_r \Psi\partial_\theta \Psi + \frac1 r \overline{U}_{2}' \{\frac{\gamma+1}{2}(\partial_r\Psi)^2+ \frac{\gamma-1}{2r^2} (\partial_{\theta}\Psi)^2\}\\\nonumber
&&+ \frac{\gamma-1}{2r\tan\theta} \overline{U}_{2} |\nabla \Psi|^2 -\frac{1}{r}(c^2(H)-c^2(\overline{\rho}))(2\partial_r\Psi+\frac{1}{r\tan\theta} \partial_{\theta} \Psi)\\\nonumber
&&+\frac{\gamma-1}{r} \overline{U}_{1} |\nabla \Psi|^2+ \frac{1}{r^3}(2\partial_{\theta}\overline{\Phi}+ \partial_{\theta}\Psi)\partial_{r}\Psi-\frac{\overline{U}_{1}}{r^3}(\partial_{\theta}\Psi)^2.
\end{eqnarray}

\subsection{The existence of a multiplier}
For small values of $\Psi$, the equation in \eqref{s100} exhibits a nonlinear mixed-type second-order equation of Tricomi type. This can be understood through the following analysis. Initially, we examine the linearized mixed-type second-order differential equation in $\Omega$, which is derived by setting $\Psi\equiv 0$ in the coefficients $A_{11}$, $A_{12}$, and $A_{22}$:
\be\label{b1}
&&\ol{\mathcal{L}}\Psi:=\overline{A}_{11}(\theta) \p_r^2 \Psi+ \frac{2 \overline{A}_{12}(\theta)}{r} \p_{r\theta}^2 \Psi + \frac{\overline{A}_{22}(\theta)}{r^2}\p_{\theta}^2\Psi\\\no
&&\quad \quad \quad\quad\quad + \frac{e_1(\theta)}{r} \p_r \Psi+ \frac{e_2(\theta)}{r^2}\p_{\theta} \Psi= F(r,\theta),
\ee
where
\be\label{ba11}\begin{cases}
\overline{A}_{11}(\theta)= c^2(\overline{\rho})- \overline{U}_{1}^2>0,\ \ \overline{A}_{12}(\theta)= -\overline{U}_{1}\overline{U}_{2}<0,\\
\overline{A}_{22}(\theta)= c^2(\overline{\rho})- \overline{U}_{2}^2>0.
\end{cases}\ee
To examine the types of equation \eqref{b1}, we can perform a change of variables to eliminate the coefficient of the mixed second-order derivative:
\be\label{b2}
y_1=\frac{f(\theta)}{r},\ \  y_2 =\theta,
\ee
where $f(\theta)= \text{e}^{\int_{\theta_c}^{\theta}\frac{\overline{A}_{12}}{\overline{A}_{22}}(\tau) d\tau}$ satisfying $\frac{f'(\theta)}{f(\theta)}=\frac{-\overline{M}_{1}\overline{M}_{2}}{1-\overline{M}_{2}^2}(\theta)$. Then $r=\frac{f(y_2)}{y_1}$ and
\be\no
&&\p_r=-\frac{y_1^2}{f(y_2)} \p_{y_1},\ \ \p_{\theta}= \frac{y_1 f'}{f}\p_{y_1} + \p_{y_2},\\\no
&&\p_r^2= \frac{y_1^4}{f^2}\p_{y_1}^2 + \frac{2 y_1^3}{f^2}\p_{y_2},\\\no
&&\p_{r\theta}^2= -\frac{y_1^2 f'}{f^2}\p_{y_1}-\frac{y_1^3 f'}{f^2}\p_{y_1}^2- \frac{y_1^2}{f(y_2)}\p_{y_1y_2}^2,\\\no
&&\p_{\theta}^2=\p_{y_2}^2+ \frac{y_1^2 (f')^2}{f^2}\p_{y_1}^2+ \frac{2 y_1 f'}{f}\p_{y_1 y_2}^2 + y_1
\big[\big(\frac{f'}{f}\big)'+\big(\frac{f'}{f}\big)^2\big]\p_{y_1}.
\ee
By defining $\psi(y_1,y_2)= \Psi\big(\frac{f(y_2)}{y_1},y_2\big)$, the equation \eqref{b1} can be rewritten as:
\be\label{b3}
\frac{y_1^2 \overline{A}_{22}}{f^2}\big(y_1^2 \overline{k}_{11}(y_2)\p_{y_1}^2+ \p_{y_2}^2+y_1 \overline{k}_{1}(y_2) \p_{y_1}+ \overline{k}_{2}(y_2)\p_{y_2} \big)\psi=\widetilde{F}(y_1,y_2),
\ee
where
\begin{eqnarray}\label{k11}\begin{cases}
\overline{k}_{11}(y_2)=\frac{1-\overline{M}_{1}^2-\overline{M}_{2}^2}{(1-\overline{M}_{2}^2)^2},\\
\overline{k}_{1}(y_2)= \frac{2(1-\overline{M}_{1}^2)}{1-\overline{M}_{2}^2}-\frac{\overline{M}_{1}^2 \overline{M}_{2}^2}{(1-\overline{M}_{2}^2)^2}-\frac{d}{dy_2}\big(\frac{\overline{M}_{1}\overline{M}_{2}}{1-\overline{M}_{2}^2}\big)-
\frac{\overline{M}_{1}\overline{M}_{2}}{1-\overline{M}_{2}^2}\frac{e_2}{\overline{A}_{22}}-\frac{e_1}{\overline{A}_{22}},\\
\overline{k}_{2}(y_2)=\frac{e_2(y_2)}{\overline{A}_{22}(y_2)}, \ \ \ \widetilde{F}(y_1,y_2)=F(\frac{1}{y_1} f(y_2), y_2).
\end{cases}\end{eqnarray}
Since $\overline{k}_{11}(y_2)$ changes sign when crossing the sonic surface $\theta=\theta_{so}$, the equation \eqref{b3} can be classified as a linear mixed-type second-order equation of Tricomi type. It is important to note that the transformation $(r,\theta)\to (y_1,y_2)$ is nonlinear, and as a result, the rectangle $\Omega$ is transformed into a distorted domain denoted as
$$D:=\big\{(y_1,y_2): \frac{1}{r_1} f(\theta)<y_1<\frac{1}{r_0} f(\theta), \theta_-<\theta<\theta_+\big\}.$$

To avoid the challenges arising from the irregular boundary of $D$, we temporarily assume that the function $\psi$ has a compact support contained within $D$. Our next objective is to identify a suitable multiplier for \eqref{b3} that allows us to establish an $H^1(\Omega)$ energy estimate for any smooth solution $\psi$ of \eqref{b3}. Taking inspiration from the multiplier used for \eqref{b3}, we can find a corresponding multiplier for \eqref{b1} in the $(r,\theta)$ coordinates and effectively handle the difficulties caused by the boundary term.

A crucial observation in identifying a suitable multiplier is the set of identities satisfied by the background solution.
\begin{proposition}\label{key}
{\it Let $(\overline{U}_{1}, \overline{U}_{2})$ denote the unique solution as stated in Proposition \ref{transonic-bg}. In this case, the following relationship holds: for all $ y_2\in [\theta_-,\theta_+]$,
 \begin{eqnarray}\label{key1}
 &&\overline{k}_{1}(y_2)\equiv 0,\\\no
 &&\overline{k}_{11}'+ 2\overline{k}_{11} \overline{k}_{2}= \frac{\overline{M}_{1}\overline{M}_{2}}{(1-\overline{M}_{2}^2)^3} (2+(\gamma-1)|\ol{\bf M}|^2)\\\label{key2}
 &&\quad\quad+ \frac{\cot y_2}{(1-\overline{M}_{2}^2)^3} \left(2(1-\overline{M}_{1}^2)+(\gamma-1)\overline{M}_{2}^2|\ol{\bf M}|^2\right)>0.
\end{eqnarray}
}\end{proposition}

\begin{proof}
To prove \eqref{key1}, we can deduce from equations \eqref{s4} and \eqref{s5} that:
\be\no
(\overline{M}_{1}\overline{M}_{2})'(y_2)= \overline{M}_{2}^2-\overline{M}_{1}^2-\frac{(\gamma-1)\overline{M}_{1} \overline{M}_{2}^2 (\overline{M}_{1}\sin y_2 + \overline{M}_{2}\cos y_2)}{2(1-\overline{M}_{2}^2)\sin y_2}
\ee
and
\be\no
&&\overline{k}_{1}=\frac{2(1-\overline{M}_{1}^2)(1-\overline{M}_{2}^2)-\overline{M}_{1}^2\overline{M}_{2}^2}{(1-\overline{M}_{2}^2)^2}-\frac{(\overline{M}_{1}\overline{M}_{2})'}{1-\overline{M}_{2}^2}
-\frac{\overline{M}_{1}\overline{M}_{2}\frac{d}{dy_2}\overline{M}_{2}^2}{(1-\overline{M}_{2}^2)^2}\\\no
&&\quad-\frac{1}{1-\overline{M}_{2}^2}\big(2-\overline{M}_{2}^2-(\gamma-1)\overline{M}_{1}^2+ \frac{(\gamma-1)\overline{M}_{1}^2}{1-\overline{M}_{2}^2}+\frac{(\gamma-1)\overline{M}_{1}\overline{M}_{2}^3\cos y_2}{(1-\overline{M}_{2}^2)\sin y_2}\big)\\\no
&&-\frac{\overline{M}_{1}\overline{M}_{2}}{(1-\overline{M}_{2}^2)^2}\big((3-\gamma)\overline{M}_{1}\overline{M}_{2}+ \frac{(\gamma+1)\overline{M}_{2}(\overline{M}_{1}\sin y_2+ \overline{M}_{2} \cos y_2)}{(1-\overline{M}_{2}^2)\sin y_2}\\\no
&&\quad\ -(\gamma-1)\overline{M}_{2}^2\cot y_2 + \cot y_2\big)\\\no
&&=\frac1{(1-\overline{M}_{2}^2)^2}\bigg\{2-2\overline{M}_{1}^2-2 \overline{M}_{2}^2+ \overline{M}_{1}^2 \overline{M}_{2}^2- (1-\overline{M}_{2}^2)(\overline{M}_{2}^2-\overline{M}_{1}^2)\\\no
&&\quad\ +\frac{\overline{M}_{1}(\overline{M}_{1}\sin y_2+ \overline{M}_{2}\cos y_2)}{\sin y_2}\big(1+(\gamma-1)\overline{M}_{2}^2\big)\\\no
&&\quad\  + 2 \overline{M}_{1}^2 \overline{M}_{2}^2 + \frac{ \overline{M}_{1}\overline{M}_{2}^2 (\overline{M}_{1}\sin y_2+ \overline{M}_{2}\cos y_2)}{(1-\overline{M}_{2}^2)\sin y_2}(2+(\gamma-1)\overline{M}_{2}^2)\\\no
&&\quad\ -(1-\overline{M}_{2}^2)(2-(\gamma-1)\overline{M}_{1}^2-\overline{M}_{2}^2)-(\gamma-1)\overline{M}_{1}^2-(\gamma-1)\overline{M}_{1}\overline{M}_{2}^3 \cot y_2\\\no
&&\quad\  -(3-\gamma)\overline{M}_{1}^2 \overline{M}_{2}^2+ \overline{M}_{1}\overline{M}_{2}((\gamma-1)\overline{M}_{2}^2-1)\cot y_2\\\no
&&\quad\ -\frac{(\gamma+1)\overline{M}_{1}\overline{M}_{2}^2(\overline{M}_{1}\sin y_2+ \overline{M}_{2}\cos y_2)}{(1-\overline{M}_{2}^2)\sin y_2}\bigg\}\equiv 0.
\ee
To complete the verification of \eqref{key2}, we need to demonstrate the following:
\be\no
&&\overline{k}_{11}'(y_2)=-\frac{1}{(1-\overline{M}_{2}^2)^2}\frac{d}{dy_2} |\ol{\bf M}|^2+ 2\frac{1-\overline{M}_{1}^2-\overline{M}_{2}^2}{(1-\overline{M}_{2}^2)^3}\frac{d}{dy_2}\overline{M}_{2}^2\\\no
&&=\frac{\overline{M}_{2}(\overline{M}_{1}\sin y_2+ \overline{M}_{2}\cos y_2)}{(1-\overline{M}_{2}^2)^3\sin y_2}(2+(\gamma-1)(\overline{M}_{1}^2+\overline{M}_{2}^2))\\\no
&&-\frac{2\overline{M}_{2}(1-\overline{M}_{1}^2-\overline{M}_{2}^2)}{(1-\overline{M}_{2}^2)^3}\big(2\overline{M}_{1}+\frac{(\overline{M}_{1}\sin y_2+ \overline{M}_{2}\cos y_2)}{(1-\overline{M}_{2}^2)\sin y_2}(2+(\gamma-1)\overline{M}_{2}^2)\big)\\\no
&&=-\frac{4\overline{M}_{1}\overline{M}_{2}(1-\overline{M}_{1}^2-\overline{M}_{2}^2)}{(1-\overline{M}_{2}^2)^3}+\frac{\overline{M}_{2}(\overline{M}_{1}\sin y_2+ \overline{M}_{2}\cos y_2)}{(1-\overline{M}_{2}^2)^4\sin y_2}\bigg((\gamma+3)\overline{M}_{1}^2\\\no
&&\q\q\q+ (3-\gamma) \overline{M}_{2}^2+(\gamma-1)\overline{M}_{2}^2(\overline{M}_{1}^2+\overline{M}_{2}^2)-2\bigg)
\ee
and
\be\no
&&\overline{k}_{11}'+ 2\overline{k}_{11} \overline{k}_{2}\\\no
&&= \frac{2(1-\overline{M}_{1}^2-\overline{M}_{2}^2)}{(1-\overline{M}_{2}^2)^3}((1-\gamma)\overline{M}_{1}\overline{M}_{2}+ (1-(\gamma-1)\overline{M}_{2}^2)\cot y_2)\\\no
&&\quad + \frac{\overline{M}_{2}(\overline{M}_{1}\sin y_2+ \overline{M}_{2}\cos y_2)}{(1-\overline{M}_{2}^2)^4\sin y_2}\bigg((\gamma+3)\overline{M}_{1}^2+ (3-\gamma) \overline{M}_{2}^2\\\no
&&\quad+(\gamma-1)\overline{M}_{2}^2(\overline{M}_{1}^2+\overline{M}_{2}^2)-2+2(\gamma+1)(1-\overline{M}_{1}^2-\overline{M}_{2}^2)\bigg)\\\no
&&:= \frac{\overline{M}_{1}\overline{M}_{2}}{(1-\overline{M}_{2}^2)^4} I_1 + \frac{\cot y_2}{(1-\overline{M}_{2}^2)^4} I_2,
\ee
where
\be\no
&&I_1=2(1-\gamma)(1-|\ol{\bf M}|^2)(1-\overline{M}_{2}^2)+2\gamma+(1-\gamma)\overline{M}_{1}^2+(1-3\gamma) \overline{M}_{2}^2\\\no
&&\quad\quad\quad + (\gamma-1) \overline{M}_{2}^2|\ol{\bf M}|^2= [2+(\gamma-1)|\ol{\bf M}|^2](1-\overline{M}_{2}^2),\\\no
&&I_2=2(1-|\ol{\bf M}|^2)(1-(\gamma-1)\overline{M}_{2}^2)(1-\overline{M}_{2}^2)\\\no
&&\q\q+\overline{M}_{2}^2\big(2\gamma+(1-\gamma) \overline{M}_{1}^2+(1-3\gamma) \overline{M}_{2}^2+(\gamma-1)\overline{M}_{2}^2|\ol{\bf M}|^2\big)\\\no
&&=[2(1-\overline{M}_{1}^2)+(\gamma-1)\overline{M}_{2}^2|\ol{\bf M}|^2](1-\overline{M}_{2}^2).
\ee
Therefore, equation \eqref{key2} has been successfully proven.

\end{proof}

Based on \eqref{key1}, we can simplify equation \eqref{b3} to the following form:
\be\label{b5}
y_1^2 \overline{k}_{11}(y_2) \p_{y_1}^2 \psi + \p_{y_2}^2 \psi + \overline{k}_{2}(y_2) \p_{y_2} \psi= G(y_1,y_2),
\ee
where $G(y_1,y_2)=\frac{f^2(y_2)}{y_1^2 \overline{A}_{22}(y_2)} \widetilde{F}(y_1,y_2)$. Moreover, we can establish the existence of multipliers in the following proposition:

\begin{proposition}
{\it There exists a multiplier $\mathcal{M}\psi= y_1^{-1} d_1(y_2) \p_{y_1} \psi + y_1^{-2} d_2(y_2) \p_{y_2}\psi$ with smooth functions $d_1(y_2)$ and $d_2(y_2)$ defined on $[\theta_-,\theta_+]$ such that for any smooth function $\psi$ with compact support in $D$, there holds
\begin{eqnarray}\label{esti1}
&&\iint_{D} (y_1^2 \overline{k}_{11} \p_{y_1}^2 \psi + \p_{y_2}^2 \psi + \overline{k}_{2}(y_2) \p_{y_2} \psi) \mathcal{M}\psi dy_1 dy_2\\\no
&&\geq 2\iint_{D} (\p_{y_1}\psi)^2+ (\frac{1}{y_1} \p_{y_2}\psi)^2 dy_1 dy_2.
\end{eqnarray}
}\end{proposition}

\begin{proof}

Take the multiplier $\mathcal{M}\psi= y_1^{\mu} d_1(y_2) \p_{y_1} \psi + y_1^{\mu-1} d_2(y_2) \p_{y_2}\psi$, where $\mu$ is a constant to be determined, $d_1, d_2$ are smooth functions of $y_2$ to be selected later. Assume that $\psi$ has compact support in $D$, then some computations yield
\be\no
\iint_D G \mathcal{M}\psi= \iint_D y_1^{\mu-1}\big(\widetilde{K}_1(y_1\p_{y_1}\psi)^2+ \widetilde{K}_2 y_1\p_{y_1}\psi \p_{y_2}\psi+ \widetilde{K}_3(\p_{y_2}\psi)^2\big)dy_1 dy_2,
\ee
where
\be\no\begin{cases}
\widetilde{K}_1(y_2)=\frac{1}{2}[ (\overline{k}_{11}d_2)'-(\mu+2) \overline{k}_{11}d_1],\\
\widetilde{K}_2(y_2)= \overline{k}_{2}d_1-d_1'-(\mu+1) \overline{k}_{11}d_2,\\
\widetilde{K}_3(y_2)= \overline{k}_{2} d_2-\frac{1}{2}d_2'+ \frac{1}{2}\mu d_1.
\end{cases}\ee
Let $\mu=-1$ and take
\be\no
d_1(y_2)= d_1(\theta_{so}) e^{\int_{\theta_{so}}^{y_2} \overline{k}_{2}(\tau) d\tau},
\ee
with any given constant $d_1(\theta_{so})$, then
\be\no
\widetilde{K}_2(y_2)= \overline{k}_{2}d_1-d_1'-(\mu+1) k_{11}d_2\equiv 0,\ \ \forall y_2\in [\theta_-,\theta_+].
\ee
Taking
\be\no
d_2(y_2)=d_2(\theta_{so}) e^{\int_{\theta_{so}}^{y_2} 2\overline{k}_{2}(\tau)d\tau}- 2 d_* \int_{\theta_{so}}^{y_2} e^{2\int_{s}^{y_2} \overline{k}_{2}(\tau)d\tau} ds,
\ee
so that $d_2'(y_2)= 2 \overline{k}_{2} d_2(y_2)-2 d_*$ for all $y_2\in [\theta_-,\theta_+]$ where $d_*:=\max\{\|d_1\|_{L^\infty},4\}$ is a positive constant such that
\be\no
\widetilde{K}_3(y_2)=\overline{k}_{2} d_2-\frac12 d_2'-\frac12 d_1=d_*-\frac 12 d_1\geq \frac 12 d_*\geq 2.
\ee
Note that
\be\no
(\overline{k}_{11}d_2)'(y_2)= \overline{k}_{11}'d_2+ \overline{k}_{11} d_2'=(\overline{k}_{11}'+ 2 \overline{k}_{11} \overline{k}_{2}) d_2- 2\overline{k}_{11} d_*.
\ee
By utilizing \eqref{key2}, for any given constant $d_1(\theta_{so})$, it is possible to select a sufficiently large positive constant $d_2(\theta_{so})$ such that
\be\no
\widetilde{K}_1(y_2)=\frac{1}{2}[(\overline{k}_{11}d_2)'- \overline{k}_{11}d_1]=\frac{1}{2}[(\overline{k}_{11}'+ 2 \overline{k}_{11} \overline{k}_{2}) d_2- \overline{k}_{11} (2d_*+ d_1)]\geq 2.
\ee
Therefore, the proof of \eqref{esti1} is completed.

\end{proof}

\subsection{The admissible boundary conditions and energy estimates to the equation \eqref{b1}}

Returning to the original coordinate system of $(r, \theta)$, the corresponding multiplier is expected to be
\begin{eqnarray}\label{m3}
\mathcal{M}_1 \Psi= r^{\mu_1} l_1(\theta)\p_r \Psi + r^{\mu_1-1} l_2(\theta) \p_{\theta} \Psi, \ \ (r,\theta)\in (r_0,r_1)\times (\theta_-,\theta_+)
\end{eqnarray}
where $\mu_1=3$ and
\be\label{m11}
l_1(\theta)= \frac{1}{f(\theta) \overline{A}_{22}(\theta)}\left(\frac{f'(\theta)}{f(\theta)} d_2(\theta)-d_1(\theta)\right),\ \ \ l_2(\theta)= \frac{d_2(\theta)}{f(\theta) \overline{A}_{22}(\theta)}.
\ee
By fixing $d_1(\theta_{so})>0$ and selecting $d_2(\theta_{so})>0$ sufficiently large such that $d_1(\theta)>0$ and $d_2(\theta)>0$, and noting that  $\frac{f'(\theta)}{f(\theta)}=\frac{-\overline{M}_{1}(\theta)\overline{M}_{2}(\theta)}{(1-\overline{M}_{2}^2)(\theta)}<0$, we can conclude that
\be\label{l_1_l_2_sign}
l_1(\theta)<0,\text{ and } l_2(\theta)>0.
\ee

To establish a well-posed boundary value problem for \eqref{b1}, it is necessary to determine a set of admissible boundary conditions. For any smooth function $\Psi$ defined on $\Omega$, the following equality holds:
\begin{eqnarray}\nonumber
&&\iint_{\Omega} F\mathcal{M}_1\Psi dr d\theta= \iint_{\Omega} \ol{\mathcal{L}}\Psi\cdot \mathcal{M}_1\Psi dr d\theta\\\nonumber
&&=\int_{\theta_-}^{\theta_+} r^{\mu_1}\big(\overline{A}_{11}\partial_r\Psi(\frac{l_1}2\partial_r\Psi+\frac{l_2}{r}\partial_{\theta}\Psi+ (\overline{A}_{12}l_2-\frac12 \overline{A}_{22} l_1)(\frac 1r \partial_{\theta}\Psi)^2\big)\bigg|_{r=r_0}^{r_1} d\theta\\\nonumber
&&+ \int_{r_0}^{r_1} r^{\mu_1-1}\big((\overline{A}_{12}l_1-\frac12 \overline{A}_{11}l_2)(\partial_r\Psi)^2+ \frac{\overline{A}_{22}}{r}\p_{\theta}\Psi(l_1 \partial_r\Psi+ \frac{l_2}{2r}\partial_{\theta}\Psi\big)\bigg|_{\theta=\theta_-}^{\theta_+} dr\\\nonumber
&&\quad + \iint_{\Omega} r^{\mu_1-1} \big(K_1(\partial_r \Psi)^2 + K_2 \partial_r \Psi\frac{1}{r}\partial_{\theta}\Psi+K_3(\frac{1}{r}\partial_{\theta}\Psi)^2\big) dr d\theta\\\nonumber
&&=\frac12\int_{\theta_-}^{\theta_+} r^{\mu_1}\big(\overline{A}_{11} l_1 (\partial_r \Psi+ \frac{l_2}{l_1}\frac{1}{r}\partial_{\theta}\Psi)^2+ K_4 (\frac{1}{r}\partial_{\theta}\Psi)^2\big)\bigg|_{r=r_0}^{r_1} d\theta \\\nonumber
&&\quad+ \frac{1}{2}\int_{r_0}^{r_1} \frac{l_1}{l_2}r^{\mu_1-1}\big(K_4(\partial_r\Psi)^2+ \overline{A}_{22}l_1(\partial_r \Psi+ \frac{l_2}{l_1}\frac{1}{r}\partial_{\theta}\Psi)^2\big)\bigg|_{\theta=\theta_-}^{\theta_+} dr\\\label{m12}
&&\quad+ \iint_{\Omega} r^{\mu_1-1} \big(K_1(\partial_r \Psi)^2 + K_2 \partial_r \Psi\frac{1}{r}\partial_{\theta}\Psi+K_3(\frac{1}{r}\partial_{\theta}\Psi)^2\big) dr d\theta,
\end{eqnarray}
where
\begin{eqnarray}\label{m121}\begin{cases}
K_1(\theta)= \frac12 \frac{d}{d\theta}(\overline{A}_{11}l_2)-\frac{1}{2} \mu_1 \overline{A}_{11} l_1 -\frac{d}{d\theta}(\overline{A}_{12} l_1) + e_1 l_1,\\
K_2(\theta)= e_1 l_2 +e_2l_1 -(\overline{A}_{22}l_1)' -(\mu_1-1)l_2\overline{A}_{11},\\
K_3(\theta)= \frac{\mu_1-2}{2} \overline{A}_{22} l_1 + e_2 l_2- (\mu_1-2) \overline{A}_{12} l_2 -\frac12 \frac{d}{d\theta} (\overline{A}_{22} l_2),\\
K_4(\theta)= 2 \overline{A}_{12} l_2 - \frac{\overline{A}_{11}l_2^2}{l_1}- \overline{A}_{22}l_1=-\frac{1}{l_1}(\overline{A}_{11}l_2^2-2\overline{A}_{12}l_1 l_2+ \overline{A}_{22}l_1^2).
\end{cases}\end{eqnarray}
In order to obtain an energy estimate, it is necessary to prove that for $\mu_1=3$ and functions $l_1$ and $l_2$ chosen as specified in \eqref{m11}, the following conditions are satisfied:

\begin{proposition}\label{prop_3_4}
For all $\theta\in [\theta_-,\theta_+]$,
\begin{eqnarray}\label{m13}
K_1(\theta)>0,\ \  K_3(\theta)>0, \ \ 4K_1(\theta) K_3(\theta) - K_2^2(\theta)>0;
\end{eqnarray}
and there exists a small positive constant $\sigma_0$ such that, for all $\theta\in [\theta_{so}-\sigma_0, \theta_+]$,
\begin{eqnarray}\label{m14}
K_4(\theta)>0.
\end{eqnarray}
\end{proposition}

\begin{proof}
After evaluating the expressions, we obtain
\be\nonumber
&&K_3(\theta)
= \frac{\mu_1-2}{2} \frac{f'(\theta)}{f^2(\theta)} d_2- \frac{\mu_1-2}{2 f(\theta)} d_1 + \frac{1}{f(\theta)}(\overline{k}_{2} d_2-\frac12 d_2')+ \frac{f'}{2 f^2} d_2\\\no
&&=\frac{1}{f(\theta)}(\overline{k}_{2} d_2-\frac12 d_2'-\frac 12 d_1)=\frac{\widetilde{K}_3(\theta)}{f(\theta)}>0,
\\\no
&&K_2(\theta)= \frac{d_2}{f(\theta)}\frac{e_1}{\overline{A}_{22}} + \frac{1}{f} \overline{k}_{2} \big(\frac{f'}{f} d_2-d_1\big)-\big(\frac{1}{f}\big(\frac{f'}{f} d_2-d_1\big)\big)'-\frac{2 \overline{A}_{11} d_2}{f \overline{A}_{22}}\\\no
&&= \frac{1}{f}\big[d_2(\frac{e_1}{\overline{A}_{22}}+\frac{f'}{f} \overline{k}_{2})-f\big(\frac{f'}{f^2} d_2\big)'- \frac{2\overline{A}_{11}d_2}{\overline{A}_{22}}\big]-\frac{\overline{k}_{2} d_1}{f}+ \big(\frac{d_1}{f}\big)'\\\no
&&= \frac{1}{f} \big(d_2\big[\frac{e_1}{\overline{A}_{22}}+\frac{\overline{k}_{2} f'}{f} -\big(\frac{f'}{f}\big)'+\big(\frac{f'}{f}\big)^2-\frac{2 \overline{A}_{11}}{\overline{A}_{22}}\big]-\frac{d_2' f'}{f}\big)+\frac{d_1'-k_2 d_1}{f}-\frac{f'd_1}{f^2}\\\no
&&=\frac{2f'}{f^2}(\overline{k}_2 d_2-\frac12d_2'-\frac12d_1)-\frac{\overline{k}_{2}d_1-d_1'}{f}=\frac{2f'\widetilde{K}_3}{f^2}-\frac{\widetilde{K}_2}{f}=\frac{2f'\widetilde{K}_3}{f^2}<0,
\ee
and
\be\no
&&K_1(\theta)
=\frac12 \frac{d}{d\theta}\big(\frac1{f}\big(\frac{\overline{A}_{11}}{\overline{A}_{22}} -\frac{\overline{A}_{12}^2}{\overline{A}_{22}^2}\big)d_2\big)-\frac{1}{2}\frac{d}{d\theta} \big(\frac{1}{f}\big(\frac{f'}{f}\big)^2 d_2\big)- \frac{3 \overline{A}_{11}}{2 f \overline{A}_{22}}\frac{f'}{f} d_2\\\no
&& \quad+ \frac{e_1f'}{f^2 \overline{A}_{22}}d_2 + \frac{3 \overline{A}_{11}}{2 f \overline{A}_{22}} d_1+\frac{d}{d\theta}\big(\frac{f' d_1}{f^2}\big)-\frac{e_1d_1}{f \overline{A}_{22}}\\\no
&&=\frac12 \big(\frac{\overline{k}_{11} d_2}{f}\big)'-\frac{f'}{f^2}\big(\frac{f'}{f}\big)' d_2 + \frac12 \frac{f'}{f^2}\big(\frac{f'}{f}\big)^2 d_2-\frac{1}{2 f} \big(\frac{f'}{f}\big)^2 d_2' + \frac{f' e_1 d_2}{f^2 \overline{A}_{22}}\\\no
&&\quad-\frac{3f' \overline{A}_{11} d_2}{2 f^2 \overline{A}_{22}}+ \frac{f'}{f^2} d_1-\frac{1}{f}\big(\frac{\overline{A}_{11}}{2\overline{A}_{22}}+\overline{k}_{2}\frac{f'}{f}\big)d_1\\\no
&&=\frac{1}{2f} (\overline{k}_{11} d_2)'+ \frac{f'}{f^2}\big[-\frac12 \overline{k}_{11}-\big(\frac{f'}{f}\big)'+\frac12 \big(\frac{f'}{f}\big)^2-\frac{3 \overline{A}_{11}}{2 \overline{A}_{22}}+\frac{e_1}{\overline{A}_{22}}\big] d_2 \\\no
&& \quad -\frac{1}{2f}\big(\frac{f'}{f}\big)^2 d_2' + \frac{f'}{f^2} d_1-\frac{1}{f}\big(\frac{\overline{A}_{11}}{2\overline{A}_{22}}+\overline{k}_{2}\frac{f'}{f}\big)d_1\\\no
&&= \frac{1}{2f} \big((\overline{k}_{11} d_2)'-\overline{k}_{11}d_1\big)+ \frac{1}{f}\big(\frac{f'}{f}\big)^2 (\overline{k}_{2} d_2-\frac12d_2'-\frac12 d_1)-\frac{f'}{f^2}(\overline{k}_{2}d_1-d_1')\\\no
&&= \frac{\widetilde{K}_1}{f}+ \frac{1}{f}\big(\frac{f'}{f}\big)^2 \widetilde{K}_3-\frac{f'}{f^2} \widetilde{K}_2=\frac{\widetilde{K}_1}{f}+ \frac{1}{f}\big(\frac{f'}{f}\big)^2 \widetilde{K}_3>0,
\\\no
&&K_4(\theta)=-\frac{1}{l_1(\theta)}(\overline{A}_{11}l_2^2-2 \overline{A}_{12} l_1 l_2 +\overline{A}_{22} l_1^2)\\\no
&&=-\frac{1}{\overline{A}_{22} f' d_2-f d_1}\big[\overline{A}_{11}d_2^2-2 \overline{A}_{12}d_2\big(\frac{\overline{A}_{12}}{\overline{A}_{22}}d_2-d_1\big)+\overline{A}_{22}\big(\frac{\overline{A}_{12}}{\overline{A}_{22}} d_2-d_1\big)^2\big]\\\no
&&=-\frac{1}{f'(y_2) d_2-f(y_2)d_1}\big(d_1^2+ d_2^2\frac{1-\overline{M}_{1}^2-\overline{M}_{2}^2}{(1-\overline{M}_{2}^2)^2}\big).
\ee

Furthermore, we observe that
\begin{equation}\no
K_4(\theta)>0,\quad \text{for any }\theta \in [\theta_{so}, \theta_+],
\end{equation}
and, by continuity, $K_4(\theta)>0$ for any $\theta\in [\theta_{so}-\sigma_0, \theta_+]$ with a small positive constant $\sigma_0$.
\end{proof}

As a result of \eqref{m12} and Proposition \ref{prop_3_4}, we are motivated to consider the following problem in the domain $\Omega_1=(r_0,r_1)\times (\theta_{so}-\sigma_0, \theta_+)$ with a class of admissible boundary conditions:
\begin{eqnarray}\label{m19}\begin{cases}
\ol{\mathcal{L}}\Psi:=\overline{A}_{11} \partial_r^2 \Psi + \frac{2 \overline{A}_{12}}{r} \partial_{r\theta}^2 \Psi +\frac{\overline{A}_{22}}{r^2} \partial_{\theta}^2 \Psi + \frac{e_1(\theta)}{r}\partial_r \Psi + \frac{e_2(\theta)}{r^2} \partial_{\theta} \Psi= F(r,\theta),\\
\frac{1}{r_0}\partial_{\theta}\Psi(r_0,\theta)= g_0(\theta),\\
\partial_r \Psi(r_1,\theta) + \frac{l_2(\theta)}{l_1(\theta)} \frac{1}{r_1}\partial_{\theta} \Psi(r_1,\theta)= g_1(\theta), \\
\partial_r \Psi(r,\theta_{so}-\sigma_0) + \frac{l_2(\theta_{so}-\sigma_0)}{l_1(\theta_{so}-\sigma_0)} \frac{1}{r}\partial_{\theta} \Psi(r,\theta_{so}-\sigma_0)= h_0(r),\\
\partial_r \Psi(r,\theta_+) = h_1(r).
\end{cases}\end{eqnarray}
Moreover, based on the above analysis, we can deduce the following theorem:
\begin{theorem}
{\it For any smooth solution $\Psi$ to the boundary value problem \eqref{m19}, the following energy estimate holds.
\be\nonumber
&&\iint_{\Omega_1} r^2[(\partial_r \Psi)^2+(\frac{1}{r}\partial_{\theta}\Psi)^2] dr d\theta \\\no
&&\quad+ \int_{\theta_{so}-\sigma_0}^{\theta_+} r_0^3(\p_r\Psi+ \frac{l_2}{l_1}\frac{1}{r_0}\p_{\theta}\Psi)^2(r_0,\theta)+ r_1(\p_{\theta}\Psi)^2(r_1,\theta)d\theta\\\no
&&\quad  + \int_{r_0}^{r_1} r^2(\p_r\Psi+ \frac{l_2}{l_1}\frac{1}{r}\p_{\theta}\Psi)^2(r,\theta_+)+ r^2(\p_r\Psi)^2(r,\theta_{so}-\sigma_0)dr \\\no
&&\leq C_* \big(\iint_{\Omega_1} r^4F^2(r,\theta) dr d\theta + \int_{\theta_{so}-\sigma_0}^{\theta_+} r_1^3 g_1^2(\theta)+ r_0^3 g_0^2(\theta) d\theta \\\label{m31}
&&\quad\quad+ \int_{r_0}^{r_1} r^2 (h_0^2(r)+h_1^2(r)) dr\big),
\ee
where $C_*$ is a positive constant that depends solely on the background solution, $\theta_{so}-\sigma_0$, and $\theta_+$, but not on $r_0$ and $r_1$. Moreover, the smooth solution to \eqref{m19} is unique.
}\end{theorem}
\begin{proof}
By\eqref{m13}, there exist a positive constant $\lambda_0$ depending only on the background solution such that
\begin{align}\no
&\iint_{\Omega_1} r^2 \big(K_1(\partial_r \Psi)^2 + K_2 \partial_r \Psi\frac{1}{r}\partial_{\theta}\Psi+K_3(\frac{1}{r}\partial_{\theta}\Psi)^2\big) dr d\theta\\\nonumber
&\geq \lambda_0 \iint_{\Omega_1} r^2 \big((\partial_r \Psi)^2 +(\frac{1}{r}\partial_{\theta}\Psi)^2\big) dr d\theta
\end{align}
In addition, by \eqref{ba11}, \eqref{l_1_l_2_sign}, \eqref{m14}, there exists a positive constant $\lambda_1$ that depends solely on the background solution, $\theta_{so}-\sigma_0$, and $\theta_+$ such that $\overline{A}_{11} l_1\leq-\lambda_1<0$, $K_4\geq\lambda_1>0$, $\frac{l_1}{l_2}K_4\leq-\lambda_1<0$, $\overline{A}_{22}\frac{l_1^2}{l_2}\geq\lambda_1>0$. Thus, by using the boundary conditions in \eqref{m19}, we can derive
\begin{eqnarray}\nonumber
&&C_0\bigg(\int_{\theta_{so}-\sigma_0}^{\theta_+} r_1^3 g_1^2(\theta)+ r_0^3 g_0^2(\theta) d\theta + \int_{r_0}^{r_1} r^2 (h_0^2(r)+h_1^2(r)) dr\bigg)\\\label{m12_new}
&&+\iint_{\Omega_1} F\mathcal{M}_1\Psi dr d\theta \geq\lambda_0\iint_{\Omega_1} r^2 \big( (\partial_r \Psi)^2 +  (\frac{1}{r}\partial_{\theta}\Psi)^2\big) dr d\theta\\\no
&&\quad+ \frac12 \lambda_1\int_{\theta_{so}-\sigma_0}^{\theta_+} r_0^3 (\p_r\Psi+ \frac{l_2}{l_1}\frac{1}{r_0}\p_{\theta}\Psi)^2(r_0,\theta)+ r_1 (\p_{\theta}\Psi)^2(r_1,\theta)d\theta\\\no
&&\quad+ \frac12 \lambda_1\int_{r_0}^{r_1} r^2 (\p_r\Psi+ \frac{l_2}{l_1}\frac{1}{r}\p_{\theta}\Psi)^2(r,\theta_+)+ r^2 (\p_r\Psi)^2(r,\theta_{so}-\sigma_0)dr,
\end{eqnarray}
where $C_0$ is a positive constant that depends solely on the background solution, $\theta_{so}-\sigma_0$, and $\theta_+$. In addition,
\begin{align}\no
&\iint_{\Omega_1} F\mathcal{M}_1\Psi dr d\theta=\iint_{\Omega_1} F\Big(r^3 l_1(\theta)\p_r \Psi + r^2 l_2(\theta) \p_{\theta} \Psi\Big) dr d\theta\\\no
\leq &\  \iint_{\Omega_1} \max\{\|l_1\|_{L^\infty},\|l_2\|_{L^\infty}\}r^2F\Big(r\p_r \Psi + \p_{\theta} \Psi\Big) dr d\theta\\\label{eq_3_51}
\leq &\ \frac{\lambda_0}2\iint_{\Omega_1} r^2 \left( (\partial_r \Psi)^2 +  (\frac{1}{r}\partial_{\theta}\Psi)^2\right) dr d\theta+\frac{C}{\lambda_0}\iint_{\Omega_1} r^4F^2 dr d\theta.
\end{align}
Therefore, by substituting \eqref{eq_3_51} into \eqref{m12_new}, we can prove the validity of equation \eqref{m31}.
\end{proof}


Next, we can consider the general linearized operator
\begin{eqnarray}\label{m33}
A_{11}\p_r^2 \Psi + \frac{2 A_{12}}{r}\p_{r\theta}^2 \Psi + \frac{A_{22}}{r^2}\p_{\theta}^2 \Psi + \frac{e_1(\theta)}{r}\p_r\Psi + \frac{e_2(\theta)}{r^2}\p_{\theta}\Psi= F_1,
\end{eqnarray}
where $A_{ij}$ for $i,j=1,2$ are $C^1$ smooth functions of $(r,\theta)$ satisfying
\begin{eqnarray}\label{m34}
\|A_{ij}-\overline{A}_{ij}\|_{C^1(\overline{\Omega})}\leq \delta_0
\end{eqnarray}
for some sufficiently small positive constant $\delta_0$. Performing the same calculations as in equation \eqref{m12}, we obtain:
\begin{eqnarray}\nonumber
&&\iint_{\Omega} F\mathcal{M}_1\Psi dr d\theta\\\no
&&=\iint_{\Omega} r^{\mu_1-1} \big(\mathcal{K}_1(\partial_r \Psi)^2 + \mathcal{K}_2 \partial_r \Psi\frac{1}{r}\partial_{\theta}\Psi+\mathcal{K}_3(\frac{1}{r}\partial_{\theta}\Psi)^2\big) dr d\theta\\\nonumber
&&+ \frac{1}{2}\int_{r_0}^{r_1} \frac{l_1}{l_2}r^{\mu_1-1}\big(\mathcal{K}_4(\partial_r\Psi)^2+ A_{22}l_1(\partial_r \Psi+ \frac{l_2}{l_1}\frac{1}{r}\partial_{\theta}\Psi)^2\big)\bigg|_{\theta=\theta_-}^{\theta_+} dr\\\nonumber
&&+\frac12\int_{\theta_-}^{\theta_+} r^{\mu_1}\big(A_{11} l_1 (\partial_r \Psi+ \frac{l_2}{l_1}\frac{1}{r}\partial_{\theta}\Psi)^2+ \mathcal{K}_4 (\frac{1}{r}\partial_{\theta}\Psi)^2\big)\bigg|_{r=r_0}^{r_1} d\theta,
\end{eqnarray}
where
\begin{eqnarray}\nonumber
&&\mathcal{K}_1(r,\theta)= \frac12 \p_{\theta}(A_{11}l_2)-\frac{1}{2} \mu_1 A_{11} l_1 -\frac{1}{2}r\p_r A_{11} l_1-\p_{\theta}(A_{12} l_1) + b_1 l_1,\\\nonumber
&&\mathcal{K}_2(r,\theta)= e_1 l_2 +e_2l_1 -\p_{\theta}(A_{22}l_1) -(\mu_1-1)l_2A_{11}-r\p_r A_{11} l_2,\\\nonumber
&&\mathcal{K}_3(r,\theta)= \frac{\mu_1-2}{2} A_{22} l_1 + e_2 l_2- (\mu_1-2) A_{12} l_2 \\\nonumber
&&\quad\quad\quad-\frac12 \p_{\theta} (A_{22} l_2)+\frac12 r\p_r A_{22} l_1- r\p_r A_{12} l_2,\\\nonumber
&&\mathcal{K}_4(r,\theta)= 2 A_{12} l_2 - \frac{A_{11}l_2^2}{l_1}- A_{22}l_1=-\frac{1}{l_1}(A_{11}l_2^2-2A_{12}l_1 l_2+ A_{22}l_1^2).
\end{eqnarray}
If the parameter $\delta_0$ in equation \eqref{m34} is sufficiently small (depending on $r_1$ as well), a similar energy estimate and uniqueness result can be established for the general operator in equation \eqref{m33}. We omit the details here and leave it to the interested reader.

Now, let's consider the existence of solutions to the boundary value problem \eqref{m19}. Due to the fact that $l_2(\theta)>0$, the boundary conditions specified in \eqref{m19} are oblique. Therefore, the Galerkin's method, as employed in \cite{ku02}, cannot be directly applied to construct approximate solutions. Instead, we will resort to the theory developed in \cite{fri58} to establish the existence of a generalized solution.

Denoting $u_1=\p_r\Psi$, $u_2=\fr1r\p_\th\Psi$ and ${\bf u}=\left( \begin{array}{cccc}
u_1 \\
u_2 \end{array} \right)$, we can rewrite equation \eqref{b1} as the following first-order linear system
\begin{align}\nonumber
\mathcal{L} {\bf u}=\left( \begin{array}{cccc}
\overline{A}_{11} & 0 \\
0 & 1 \end{array} \right)\p_r{\bf u}+\left( \begin{array}{cccc}
2A_{b12} & \overline{A}_{22} \\
-1 & 0 \end{array} \right)\fr1r\p_\th{\bf u}+\left( \begin{array}{cccc}
e_1 & e_2 \\
0 & 1 \end{array} \right)\fr1r{\bf u}=\left( \begin{array}{cccc}
F \\
0 \end{array} \right).
\end{align}
The multiplier $\mathcal{M}_1\Psi$ defined in \eqref{m3} and \eqref{m11} corresponds to the symmetrizer
\begin{eqnarray}\nonumber
Z:=2r^{\mu_1}\left( \begin{array}{cccc}
l_1(\th) & \overline{A}_{11}(\th) l_2(\th)\\
l_2(\th) & 2\overline{A}_{12}l_2(\th)-\overline{A}_{22}l_1(\th) \end{array} \right),
\end{eqnarray}
which symmetrizes the differential operator $\mathcal{L}$. It should be noted that $\det Z=4r^{2\mu_1}(2\overline{A}_{12}l_1l_2-\overline{A}_{22}l_1^2-\overline{A}_{11}l_2^2=4r^{2\mu_1}K_4l_1<0$ for any $\theta\in [\theta_{so}-\sigma_0,\theta_+]$. Denote
\begin{align*}
\mathcal{K}= Z \mathcal{L}=2\mathcal{A}^{(1)}\p_r+2\mathcal{A}^{(2)}\p_\th+\mathcal{A}^{(3)},\ \ \ {\bf f}=2r^{\mu_1}F\left( \begin{array}{cccc}
l_1 \\
l_2\end{array} \right),
\end{align*}
with
\begin{align}
&\mathcal{A}^{(1)}=r^{\mu_1}\left( \begin{array}{cccc}
l_1\overline{A}_{11} & l_2\overline{A}_{11}\\\nonumber
l_2\overline{A}_{11} & 2\overline{A}_{12}l_2-\overline{A}_{22}l_1 \end{array} \right),\\\nonumber
&\mathcal{A}^{(2)}=r^{\mu_1-1}\left( \begin{array}{cccc}
2\overline{A}_{12}l_1-\overline{A}_{11}l_2& \overline{A}_{22}l_1\\\nonumber
\overline{A}_{22}l_1 & \overline{A}_{22}l_2 \end{array} \right),\\\nonumber
&\mathcal{A}^{(3)}=2r^{\mu_1-1}\left( \begin{array}{cccc}
l_1e_1 & l_1e_2+\overline{A}_{11}l_2\\\nonumber
l_2e_1 & l_2e_2+2\overline{A}_{12}l_2-\overline{A}_{22}l_1 \end{array} \right).
\end{align}

Then $\mathcal{A}^{(i)} (i=1,2)$ are symmetric. With $l_1$ and $l_2$ defined in \eqref{m11}, one can check that the symmetric part of the matrix
\begin{eqnarray}\nonumber
\mathcal{Q}\equiv \mathcal{A}^{(3)}-\p_r \mathcal{A}^{(1)}-\p_\th \mathcal{A}^{(2)}
\end{eqnarray}
is positive definite. Indeed, we have
\begin{eqnarray}\nonumber
\fr12(\mathcal{Q}+\mathcal{Q}^T)= r^{\mu_1-1}\left( \begin{array}{cccc}
  2K_1 & K_2\\
  K_2 & 2K_3\end{array} \right)\geq  C_* I,
\end{eqnarray}
for some $C_*>0$ and $\det Z\neq0$, where $K_j(\theta)$ for $j=1,2,3$ are defined as in \eqref{m121}.

It remains to show that the boundary conditions prescribed in \eqref{m19} are admissible in the sense defined by \cite{fri58} for symmetric positive linear PDE $\mathcal{K}{\bf u}={\bf f}_1$.

Denote the outer normal vector at boundary by ${\bf n}=(n_1,n_2)^\top$, then ${\bf n} (r_0,\th)=(-1,0)^\top$, ${\bf n}(r_1,\th)=(1,0)^\top$,  ${\bf n}(r,\th_0)=(0,-1)^\top$,  ${\bf n}(r,\th_1)=(0,1)^\top$. Define $\mathcal{B}=n_1\mathcal{A}^{(1)}+n_2\mathcal{A}^{(2)}$ and we should find an appropriate decomposition $\mathcal{B}=\mathcal{B}_+ +\mathcal{B}_-$ such that $\mathcal{G}=\mathcal{B}_+ - \mathcal{B}_-$ is nonnegative definite.

Since $K_4>0$, one can decompose $\mathcal{B}(r_0,\th)$, $\mathcal{B}(r_1,\th)$, $\mathcal{B}(r,\th_0)$ and $\mathcal{B}(r,\th_1)$ as follows.
\begin{align*}
  &\mathcal{B}_+(r_0,\th)=-r_0^{\mu_1}\overline{A}_{11}l_1\left( \begin{array}{cccc}
  1 & l_2/l_1\\
  l_2/l_1 & l_2^2/l_1^2 \end{array} \right),\ \ \mathcal{B}_-(r_0,\th)=-r_0^{\mu_1}\left( \begin{array}{cccc}
  0 & 0\\
  0 & K_4 \end{array} \right),\\
&\mathcal{B}_-(r_1,\th)=r_1^{\mu_1}\overline{A}_{11}l_1\left( \begin{array}{cccc}
  1 & l_2/l_1\\
  l_2/l_1 & l_2^2/l_1^2 \end{array} \right),\ \ \mathcal{B}_+(r_1,\th)=r_1^{\mu_1}\left( \begin{array}{cccc}
  0 & 0\\
  0 & K_4 \end{array} \right),\\
&\mathcal{B}_-(r,\theta_{so}-\sigma)=-r^{\mu_1-1}\overline{A}_{22}l_2\left( \begin{array}{cccc}
  l_1^2/l_2^2  & l_1/l_2\\
  l_1/l_2 & 1\end{array} \right)(\theta_{so}-\sigma),\\
&\mathcal{B}_+(r,\theta_{so}-\sigma)=-r^{\mu_1-1}\left( \begin{array}{cccc}
  K_4l_1/l_2 & 0\\
  0 & 0 \end{array} \right)(\theta_{so}-\sigma),\\
&\mathcal{B}_+(r,\theta_+)=r^{\mu_1-1}\overline{A}_{22}l_2\left( \begin{array}{cccc}
  l_1^2/l_2^2  & l_1/l_2\\
  l_1/l_2 & 1\end{array} \right)(\theta_+),\\
&\mathcal{B}_-(r,\theta_+)=r^{\mu_1-1}\left( \begin{array}{cccc}
  K_4l_1/l_2 & 0\\
  0 & 0 \end{array} \right)(\theta_+).
\end{align*}
Define $\mathcal{M}=\mathcal{G}-\mathcal{B}=-2\mathcal{B}_-$, one can verify that $\mathcal{M}u=-2\mathcal{B}_- {\bf u}=0$ on $\partial\Omega$ is an admissible boundary condition in the sense that there exists no subspace of $\mathbb{R}^2$ containing the null space of the matrix $\mathcal{M}$ properly in which ${\bf u}^\top \mathcal{B} {\bf u}\geq 0$. The boundary condition $\mathcal{M}{\bf u}=0$ exactly corresponds to the boundary conditions in \eqref{m19}, that is
\begin{eqnarray}\nonumber\begin{cases}
u_2(r_0,\theta)= g_0(\theta),\\
u_1(r_1,\theta) + \frac{l_2(\theta)}{l_1(\theta)}u_2(r_1,\theta)= g_1(\theta), \\
u_1(r,\theta_{so}-\sigma_0) + \frac{l_2(\theta_{so}-\sigma_0)}{l_1(\theta_{so}-\sigma_0)} u_2(r,\theta_{so}-\sigma_0)= h_0(r),\\
u_1(r,\theta_+) = h_1(r).
\end{cases}\end{eqnarray}
By performing an integration by parts, we obtain:
\begin{eqnarray}\nonumber
\int_{\Omega} {\bf u}\cdot \mathcal{K} {\bf u} drd\th + \int_{\partial\Omega} {\bf u}\mathcal{M}{\bf u} dS= \int_{\Omega} {\bf u}\cdot \mathcal{Q} {\bf u} drd\th + \int_{\partial\Omega} {\bf u}\mathcal{G}{\bf u} dS,
\end{eqnarray}
from which we can derive the following estimate
\begin{eqnarray}\nonumber
\|{\bf u}\|_{L^2(\Omega)}\leq C_*\|{\bf f}\|_{L^2(\Omega)},
\end{eqnarray}
for any functions ${\bf u}$ solving
\begin{eqnarray}\label{problem4}\begin{cases}
\mathcal{K} {\bf u} ={\bf f},\ \ \text{in} \ \Omega,\\
\mathcal{B} {\bf u} = {\bf 0},\ \ \text{on}\ \partial\Omega.
\end{cases}\end{eqnarray}
Therefore, based on the results established in \cite{fri58}, it can be concluded that:
\begin{theorem}\label{weak}
{\it For any ${\bf f}\in L^2(\Omega)$, there exists a weak solution ${\bf u}\in L^2(\Omega)$ to \eqref{problem4} in the sense that
\begin{eqnarray}\nonumber
\int_{\Omega} {\bf f}\cdot {\bf v} drd\th= \int_{\Omega} \mathcal{K}^* {\bf v} \cdot {\bf u} drd\th
\end{eqnarray}
holds for any ${\bf v}\in C^1(\Omega)$ satisfying $\mathcal{M}^* {\bf v}={\bf 0}$ on $\partial\Omega$. Here
$$\nonumber
\mathcal{K}^* {\bf v}=-\mathcal{A}^{(1)}\partial_x {\bf v}- \mathcal{A}^{(2)}\partial_y {\bf v}-\p_x(\mathcal{A}^{(1)}{\bf v})-\p_y(\mathcal{A}^{(2)}{\bf v})+ \mathcal{Q}^\top {\bf v},
$$
and $\mathcal{M}^*{\bf v}= \mathcal{G}^\top {\bf v}+ \mathcal{B}{\bf v}$.
}\end{theorem}

\subsection{The linearized stability analysis of three dimensional smooth irrotational transonic flows}

In this subsection, we investigate perturbations within the class of three-dimensional irrotational flows. In the case of irrotational flows, there exists a potential function $\Phi=\Phi(r,\theta,\varphi)$ such that
\be\no
U_1=\p_r \Phi,\ \ U_2=\frac{1}{r}\p_{\theta} \Phi,\ \ \ U_3=\frac{1}{r\sin\theta}\p_{\varphi} \Phi.
\ee
By substituting the expression of the density given by \eqref{density-formula} into the first equation of the system \eqref{SCEQF}, we obtain:
\begin{eqnarray}\nonumber
&&\left[c^2(H)-(\partial_r\Phi)^2\right] \partial_r^2\Phi-\frac{2}{r^2}\partial_r \Phi \partial_{\theta}\Phi \partial_{r\theta}^2\Phi+ \left[c^2(H)-(\frac1r\partial_{\theta}\Phi)^2\right] \frac{1}{r^2}\partial_{\theta}^2\Phi \\\nonumber
&&-\frac{2}{r^2\sin\theta}\p_{\theta}\Phi\p_{\varphi}\Phi\frac{1}{r^2\sin\theta}\p_{\theta\varphi}^2\Phi
+\left[c^2(H)-\left(\frac{1}{r\sin\theta}\partial_{\varphi}\Phi\right)^2\right]\frac{1}{r^2\sin^2\theta}\p_{\varphi}^2\Phi\\\nonumber
&& -\frac{2}{r^2\sin^2\theta} \p_r\Phi\p_{\varphi}\Phi \p_{r\varphi}^2\Phi + \frac2r c^2(H) \partial_r \Phi + \frac{c^2(H)}{r^2\tan\theta} \partial_{\theta}\Phi \\\nonumber
&& + \frac{1}{r^3} \partial_r \Phi (\partial_{\theta}\Phi)^2+ \frac{1}{r^3\sin^2\theta}\p_r\Phi(\p_{\varphi}\Phi)^2+\frac{\cos\theta}{r^4\sin^3\theta}\p_{\theta}(\p_{\varphi}\Phi)^2=0.
\end{eqnarray}
Also, we can denote the difference between the potential function $\Phi$ and the background solution $\overline{\Phi}$ as $\Psi=\Phi-\overline{\Phi}$. Then, the equation governing the perturbation $\Psi$ can be expressed as follows:
\begin{eqnarray}\nonumber
&&\mathcal{T}\Psi:=A_{11}(\nabla \Psi) \partial_r^2 \Psi + \frac{2 A_{12}(\nabla \Psi)}{r} \partial_{r\theta}^2 \Psi +\frac{A_{22}(\nabla \Psi)}{r^2} \partial_{\theta}^2 \Psi\\\no
&&+\frac{2 A_{13}(\nabla \Psi)}{r\sin\theta}\p_{r\varphi}^2 \Psi
\quad+\frac{2 A_{23}(\nabla \Psi)}{r^2\sin\theta} \p_{\varphi\theta}^2 \Psi + \frac{A_{33}(\nabla \Psi)}{r^2\sin^2\theta}\p_{\varphi}^2\Psi\\\no
&&+ \frac{b_1(\theta)}{r}\partial_r \Psi + \frac{b_2(\theta)}{r^2} \partial_{\theta} \Psi=F(\nabla \Psi),
\end{eqnarray}
where $A_{ij}(\nabla\Psi)$ for $i,j=1,2$ and $b_k, k=1,2$ are same as \eqref{s101} and
\begin{eqnarray}\nonumber
&&A_{13}(\nabla \Psi)= - \frac{1}{r\sin\theta}\p_{\varphi}\Psi \p_r\Psi,\ \ A_{23}(\nabla \Psi)=-\frac{1}{r^2\sin\theta}\p_{\theta}\Psi \p_{\varphi} \Psi,\\\nonumber
&& A_{33}(\nabla \Psi)= c^2(H)- \left(\frac{1}{r\sin\theta}\partial_{\varphi} \Psi\right)^2.
\end{eqnarray}
By linearizing the operator $\mathcal{T}$ at the background solution, we obtain:
\be\label{3b1}
&&\ol{\mathcal{T}}\Psi:=\overline{A}_{11}(\theta) \p_r^2 \Psi+ \frac{2 \overline{A}_{12}(\theta)}{r} \p_{r\theta}^2 \Psi + \frac{\overline{A}_{22}(\theta)}{r^2}\p_{\theta}^2\Psi \\\nonumber
&&\quad\quad+ \frac{\overline{A}_{33}(\theta)}{r^2\sin^2\theta}\p_{\varphi}^2\Psi+ \frac{e_1(\theta)}{r} \p_r \Psi+ \frac{e_2(\theta)}{r^2}\p_{\theta} \Psi= F,
\ee
where $\overline{A}_{33}(\theta)=c^2(\overline{\rho}(\theta))$. In comparison to the operator $\overline{\mathcal{L}}$ defined in \eqref{b1}, the operator $\overline{\mathcal{T}}$ defined in \eqref{3b1} incorporates an additional term $\frac{\overline{A}_{33}(\theta)}{r^2\sin^2\theta}\p_{\varphi}^2\Psi$.
Denoting the domain as $D:=(r_0,r_1)\times (\theta_{so}-\sigma_0,\theta_+)\times [0,2\pi]$, and considering the multiplier defined in equations \eqref{m3} to \eqref{m11}, the following equality can be obtained through integration by parts:
\be\no
&&\iiint_D \frac{\overline{A}_{33}}{r^2\sin^2\theta} \p_{\varphi}^2 \Psi (r^3 l_1\p_{r}\Psi+ r^2 l_2 \p_{\theta}\Psi) dr d\theta d\varphi\\\no
&&=-\frac{1}{2}\int_{\theta_{so}-\sigma_0}^{\theta_+} \int_0^{2\pi} r\frac{\overline{A}_{33}l_1}{\sin^2 \theta}(\p_{\varphi}\Psi)^2\bigg|_{r=r_0}^{r_1} d\theta d\varphi\\\no
&&-\frac12 \int_{r_0}^{r_1} \int_0^{2\pi} \frac{\overline{A}_{33}l_2}{\sin^2\theta}(\p_{\varphi}\Psi)^2\bigg|_{\theta=\theta_{so}-\sigma_0}^{\theta_+} dr d\varphi\\\no
&&+\frac12 \iiint_D \big(\overline{k}_{33}' d_2+ \overline{k}_{33}(2 \overline{k}_{2} d_2-\frac{2 d_2\cos y_2}{\sin y_2} -2 d_*- d_1\big)\frac{(\p_{\varphi}\Psi)^2}{\sin^2\theta} dr d\theta d\varphi,
\ee
where
\be\label{k33}
\overline{k}_{33}(\theta)=\frac{1}{1-\overline{M}_{2}^2(\theta)}.
\ee
In addition, we also observe the following property of the background solution:
\begin{proposition}
{\it Let $(\overline{U}_{1}, \overline{U}_{2})$ be the unique transonic solution in Proposition \ref{transonic-bg}, then for any $\theta\in [\theta_-,\theta_+]$, there holds
\begin{eqnarray}\label{key3}
\overline{k}_{33}'+ 2\overline{k}_{33} (\overline{k}_{2}-\cot \theta)=\frac{[4+(\gamma-3) \overline{M}_{2}^2]\overline{M}_{2}(\overline{M}_{1}+\overline{M}_{2}\cot \theta)}{(1-\overline{M}_{2}^2)^3}>0.
\end{eqnarray}
}\end{proposition}

\begin{proof}
By simple calculations, we have
\be\no
\overline{k}_{33}'(\theta)=-\frac{1}{(1-\overline{M}_{2}^2)^2}\big[2\overline{M}_{1}\overline{M}_{2} +\frac{\overline{M}_{2}(\overline{M}_{1}\sin\theta+\overline{M}_{2}\cos\theta)}{(1-\overline{M}_{2}^2)\sin\theta}\big(2+(\gamma-1) \overline{M}_{2}^2\big)\big]
\ee
and
\be\no
\overline{k}_{2}(\theta)=\frac{\overline{M}_{2}}{1-\overline{M}_{2}^2}\big\{[4+(\gamma-3)\overline{M}_{2}^2]\frac{(\overline{M}_{1}\sin \theta+\overline{M}_{2}\cos \theta)}{(1-\overline{M}_{2}^2)\sin \theta}-\overline{M}_{2}\cot \theta\big\}+ \cot \theta.
\ee
Then \eqref{key3} easily follows from the above two identities.
\end{proof}
Based on equation \eqref{key3}, we can choose a sufficiently large positive constant $d_2(\theta_{so})$ such that for any $\theta \in [\theta_-,\theta_+]$, the following condition holds:
\be\no
(\overline{k}_{33}'+ 2 \overline{k}_{2} \overline{k}_{33}-\frac{2 \overline{k}_{33} \cos y_2}{\sin y_2}) d_2 -2 \overline{k}_{33} d_*- \overline{k}_{33} d_1\geq C_*>0.
\ee

Motivated by the properties of the background solution, we can specify the admissible boundary conditions for equation \eqref{3b1} as follows:
\begin{eqnarray}\label{m402}\begin{cases}
\frac{1}{r_0}\partial_{\theta}\Psi(r_0,\theta,\varphi)= g_0(\theta,\varphi),\ \ \frac{1}{r_0\sin\theta} \p_{\varphi} \Psi(r_0,\theta,\psi)= q_0(\theta,\varphi),\\
\partial_r \Psi(r_1,\theta,\varphi) + \frac{l_2(\theta)}{l_1(\theta)} \frac{1}{r_1}\partial_{\theta} \Psi(r_1,\theta,\varphi)= g_1(\theta,\varphi), \\
\partial_r \Psi(r,\theta_{so}-\sigma_0,\varphi) + \frac{l_2(\theta_{so}-\sigma_0)}{l_1(\theta_{so}-\sigma_0)} \frac{1}{r}\partial_{\theta} \Psi(r,\theta_{so}-\sigma_0,\varphi)= h_0(r,\varphi),\\
\partial_r \Psi(r,\theta_+,\varphi) = h_1(r,\varphi),\ \ \frac{1}{r\sin \theta_+}\p_{\varphi} \Psi(r,\theta_+,\varphi)= q_1(r,\varphi).
\end{cases}\end{eqnarray}
Based on the analysis and discussions presented, we can now state and conclude the following theorem:
\begin{theorem}
{\it For any smooth solution $\Psi$ to the boundary value problem \eqref{3b1} and \eqref{m402}, the following energy estimate holds.
\be\nonumber
&&\iiint_{D} r^2[(\partial_r \Psi)^2+(\frac{1}{r}\partial_{\theta}\Psi)^2+(\frac{1}{r\sin\theta}\partial_{\varphi}\Psi)^2] dr d\theta d\varphi\\\nonumber
&&+ \int_{\theta_{so}-\sigma_0}^{\theta_+}\int_0^{2\pi} r_0^3(\p_r\Psi+ \frac{l_2\p_{\theta}\Psi}{r_0l_1})^2(r_0,\cdot)+ r_1[(\p_{\theta}\Psi)^2+(\frac{\p_{\varphi}\Psi}{\sin\theta})^2](r_1,\cdot)d\theta d\varphi\\\nonumber
&&+ \int_{r_0}^{r_1}\int_0^{2\pi} r^2(\p_r\Psi+ \frac{l_2\p_{\theta}\Psi}{l_1r})^2(r,\theta_+,\varphi)+[ r^2(\p_r\Psi)^2++(\frac{\p_{\varphi}\Psi}{\sin\theta})^2](r,\theta_{so}-\sigma_0,\varphi) dr d\varphi \\\no
&&\leq C_{**} \big(\iiint_{D} r^4F^2(r,\theta,\varphi) dr d\theta d\varphi + \int_{r_0}^{r_1}\int_0^{2\pi} r^2 (h_0^2+h_1^2+ q_1^2)(r,\varphi) dr d\varphi \\\no
&&\quad+ \int_{\theta_{so}-\sigma_0}^{\theta_+}\int_0^{2\pi} [r_1^3 g_1^2+ r_0^3 g_0^2+r_0^3 q_0^2](\theta,\varphi) d\theta d\varphi\big),
\ee
where $C_{**}$ is a positive constant that depends solely on the background solution, $\theta_{so}-\sigma_0$, and $\theta_+$, but not on $r_0$ and $r_1$. Moreover, the smooth solution to \eqref{3b1}-\eqref{m402} is unique.
}\end{theorem}

\begin{remark}\label{r401}
{\it One can also show the existence of generalized $L^2$ solution to \eqref{3b1} and \eqref{m402} in the sense of Friedrichs \cite{fri58} by mimicking the argument in the previous subsection.
}\end{remark}

\section{A class of smooth transonic flows with nonzero vorticity on the polar and azimuthal angles coordinates}\label{transonic2}

The goal in this section is to construct a special class of solutions to the steady non-isentropic Euler system. These solutions depend solely on the polar angle $\theta$ and the azimuthal angle $\varphi$ in spherical coordinates:
\begin{eqnarray}\label{sp00}
{\bf U}(\th,\var)=(U_1,U_2,U_3)(\th,\var),\ p=p(\th,\var),\ \rho=\rho(\th,\varphi),
\end{eqnarray}
and they are also close to the smooth self-similar irrotational transonic flows in Proposition \ref{transonic-bg}.
Within the defined class of special solutions \eqref{sp00}, the system \eqref{SCEQF} can be simplified to the following system of equations:
\begin{align}\label{bg2}
\begin{cases}
\p_\th(\rho U_2)+\frac1{\sin\th}\p_\var(\rho U_3)+2\rho U_1+\frac{\cos \th}{\sin\th}\rho U_2=0,\\
(U_2\p_\th+\frac{U_3}{\sin\th}\p_\var)U_1-(U_2^2+U_3^2)=0,\\
(U_2\p_\th+\frac{U_3}{\sin\th}\p_\var)U_2+\dfr{1}{\rho}\p_\th p+U_1 U_2-\frac{\cos \th}{\sin\th}U_3^2=0,\\
(U_2\p_\th+\frac{U_3}{\sin\th}\p_\var)U_3+\dfr{1}{\rho\sin\th}\p_\var p+U_1U_3+\frac{\cos \th}{\sin\th}U_2U_3=0,\\
(U_2\p_\th+\frac{U_3}{\sin\th}\p_\var)A=0.
\end{cases}
\end{align}
The system \eqref{bg2} is of mixed type and hyperbolic-elliptic coupled in subsonic region. To effectively decompose the hyperbolic and elliptic modes,  we employ the deformation-curl decomposition method, as introduced in \cite{wex19, w19}, for the three-dimensional steady Euler system.
Based on \eqref{bg2}, the Bernoulli's function and the entropy satisfy the following transport equations
\be\label{bg201}\begin{cases}
(U_2\p_\th+\frac{U_3}{\sin\th}\p_\var)B=0,\\
(U_2\p_\th+\frac{U_3}{\sin\th}\p_\var)A=0.
\end{cases}\ee
For the vector field ${\bf u}(x)= U_1(\theta,\var){\bf e}_r + U_2(\theta,\var){\bf e}_{\theta}+ U_3(\theta,\var){\bf e}_{\var}$, the corresponding vorticity field has a form ${\bm \omega}=\nabla \times {\bf u}= \frac{1}{r}(\omega_r {\bf e}_r + \omega_{\theta}{\bf e}_{\theta}+ \omega_{\varphi} {\bf e}_{\varphi})$, where
\begin{align}\label{bg410}
\om_r=\p_\th U_3-\dfr{1}{\sin\th}\p_\var U_2+\dfr{\cos\th}{\sin\th}U_3,\ \om_\th=\dfr{1}{\sin\th}\p_\var U_1-U_3,\ \om_\var=-\p_\th U_1+U_2.
\end{align}
Substituting \eqref{bg410} into the identity $\nabla \cdot(\nabla\times {\bf u})\equiv 0$ yields that
\begin{align}\label{bg4}
\p_\th\om_\th+\frac{1}{\sin\theta}\p_\var\om_\var+\om_r+\frac{\cos\theta}{\sin\theta}\om_\theta=0.
\end{align}
By \eqref{bg410}, the second equation in \eqref{bg2} is equivalent to
\begin{align}\label{bg3_1}
U_3\om_\th-U_2\om_\var=0,
\end{align}
and the fourth equation in \eqref{bg2} is equivalent to
\begin{align}\label{bg3_2}
U_2\om_r-U_1\om_{\theta}=\frac{-1}{\sin\theta}\p_{\varphi} B+\frac{B-\frac12 \sum_{j=1}^3 U_j^2}{A\gamma}\frac{1}{\sin\theta}\p_{\varphi} A.
\end{align}
Then
\be\label{bg5}
\omega_{\var}=\frac{U_3}{U_2} \omega_{\theta},\ \ \omega_r= \frac{U_1}{U_2}\omega_{\theta}- \frac{\p_{\varphi} B}{U_2\sin\theta}  +\frac{B-\frac12 \sum_{j=1}^3 U_j^2}{\gamma A U_2}\frac{\p_{\varphi} A}{\sin\theta}.
\ee
By substituting equation \eqref{bg5} into equation \eqref{bg4}, we can derive the transport equation satisfied by $\omega_{\theta}$:
\be\label{bg6}
&&\p_{\theta}\om_{\theta} + \frac{U_3}{U_2} \frac{1}{\sin\theta} \p_{\varphi} \omega_{\theta} + \bigg\{\frac{U_1}{U_2}+ \frac{1}{\sin\theta}\p_{\var}\big(\frac{U_3}{U_2}\big)+\frac{\cos\theta}{\sin\theta}\bigg\}\om_{\theta}\\\no
&&\quad=\frac{\p_{\varphi} B}{U_2\sin\theta}-\frac{B-\frac12 \sum_{j=1}^3 U_j^2}{\gamma A U_2}\frac{\p_{\varphi} A}{\sin\theta}.
\ee
After solving the transport equations \eqref{bg201} and \eqref{bg6} for $B$, $A$, and $\omega_{\theta}$, we can utilize \eqref{bg5} to determine $\omega_{r}$ and $\omega_{\varphi}$. Additionally, we can represent $U_2$ and $U_3$ as follows by employing \eqref{bg410}:
\be\label{bg7}\begin{cases}
U_2 =\p_{\theta} U_1+\om_{\var}=\p_{\theta} U_1+\frac{U_3}{U_2} \om_{\theta},\\
U_3 =\frac{1}{\sin\theta} \p_{\var} U_1 - \omega_{\th}.
\end{cases}\ee

It follows from the definition \eqref{Bernoulli} of the Bernoulli's function that
\begin{eqnarray}\label{bg8}
\rho= H(B,A, |{\bf U}|^2)=\big(\frac{\gamma-1}{A\gamma}\big)^{\frac{1}{\gamma-1}} \big(B-\frac12 |{\bf U}|^2\big)^{\frac{1}{\gamma-1}}.
\end{eqnarray}
By substituting \eqref{bg8} into the first equation in \eqref{bg2}, we get
\be\label{bg301}
&&(c^2(B,|{\bf U}|^2)-U_2^2)\p_\th U_2+(c^2(B,|{\bf U}|^2)-U_3^2)\frac{1}{\sin\th}\p_\var U_3\\\no
&&\quad-U_2U_3(\p_\th U_3+\frac{1}{\sin\th}\p_\var U_2)-U_1U_2\p_\th U_1-U_1U_3\frac{1}{\sin\th}\p_\var U_1\\\no
&&\quad+c^2(B,|{\bf U}|^2)(2 U_1+\frac{\cos\th}{\sin\th} U_2)=0,
\ee
where $c^2(B,|{\bf U}|^2)=(\gamma-1)(B-\frac{1}{2}\sum_{j=1}^3 U_j^2)$.
Substituting \eqref{bg7} into \eqref{bg301} yields a quasilinear second-order elliptic equation for $U_1$:
\begin{eqnarray}\nonumber
&&\big(c^2(B,|{\bf U}|^2)-(\p_{\theta} U_1+\frac{U_3}{U_2} \om_{\theta})^2\big)\p_{\theta}(\p_{\theta} U_1+\frac{U_3}{U_2}\omega_{\theta}) \\\no
&&+\big(c^2(B,|{\bf U}|^2)-(\frac{1}{\sin\theta} \p_{\var} U_1 - \omega_{\th})^2\big)\frac{1}{\sin\theta}\p_{\varphi} \big(\frac{1}{\sin\theta}\p_{\varphi}U_1-\om_{\theta}\big)\\\no
&&-(\p_{\theta} U_1+\frac{U_3}{U_2} \om_{\theta})(\frac{\p_{\varphi}U_1}{\sin\theta}  - \omega_{\th})\big(\p_{\theta}(\frac{\p_{\varphi}U_1}{\sin\theta}-\om_{\theta})+\frac{1}{\sin\theta}\p_{\varphi}(\p_{\theta} U_1+\frac{U_3}{U_2}\om_{\theta})\big)\\\no
&&-U_1(\p_{\theta} U_1+\frac{U_3}{U_2}\omega_{\theta})\p_{\theta} U_1 - U_1\big(\frac{1}{\sin\theta}\p_{\varphi} U_1-\omega_{\theta}\big)\frac{1}{\sin\theta}\p_{\varphi} U_1\\\label{bg100}
&&+ c^2(B,|{\bf U}|^2)\big(2U_1+\frac{\cos\theta}{\sin\theta}(\p_{\theta}U_1+\frac{U_3}{U_2}\omega_{\theta})\big)=0.
\end{eqnarray}

It is worth noting that the coefficient of the zeroth order term of $U_1$ in equation \eqref{bg100} (i.e. $2c^2(B,|{\bf U}|^2)$) is positive. Consequently, when applying standard boundary conditions, the equation \eqref{bg100} may not have a unique solution. Before specifying admissible boundary conditions, let us first present an equivalent reformulation of the steady Euler system \eqref{bg2} based on our previous analysis:
\begin{proposition}\label{equiv}({\bf Equivalence})
{\it Assume that $C^1$ smooth vector functions $(\rho, {\bf U}, A)$ defined on the domain $\mathbb{E}:=(\theta_-,\theta_+)\times \mathbb{T}_{2\pi}$ do not contain the vacuum (i.e. $\rho(\theta,\var)>0$ in $\mathbb{E}$) and the azimuthal velocity $U_2(\theta,\varphi)>0$ in $\mathbb{E}$, then the following two statements are equivalent:
\begin{enumerate}[(i)]
  \item $(\rho, {\bf U}, A)$ satisfy the steady Euler system \eqref{bg2} in $\mathbb{E}$;
  \item $({\bf U}, A, B)$ satisfy the equation \eqref{bg201}, \eqref{bg6}, \eqref{bg7} and \eqref{bg301} in $\mathbb{E}$ where $\om_\th:=\dfr{1}{\sin\th}\p_\var U_1-U_3$.
\end{enumerate}
}\end{proposition}

\begin{proof}
Based on our previous analysis, we have established that Statement (i) implies Statement (ii). It remains to establish the converse implication.
To proceed, let us introduce the symbol $\rho$ according to equation \eqref{bg8}, and define $\omega_\varphi$ and $\omega_r$ using equation \eqref{bg5}. By utilizing equation \eqref{bg6} and the definitions of $\omega_\varphi$ and $\omega_r$, it can be directly verified that equations \eqref{bg4}, \eqref{bg3_1}, and \eqref{bg3_2} hold. Consequently, the second and fourth equations in \eqref{bg2} are satisfied. Moreover, by combining the second and fourth equations with equation \eqref{bg201}, we can deduce that the third equation in \eqref{bg2} also holds.
Finally, the first equation in \eqref{bg2} can be derived by utilizing the definition of $\rho$ in equation \eqref{bg8}, along with equations \eqref{bg7} and \eqref{bg301}.
\end{proof}

In order to investigate the structural stability of the smooth self-similar irrotational transonic flows discussed in Proposition \ref{transonic-bg} under small perturbations of suitable boundary conditions, we need to prescribe the following boundary conditions:
\begin{eqnarray}\label{bg10}
&&B(\theta_-,\varphi)= B_0 + \epsilon B_{in}(\varphi),\\\label{bg11}
&&A(\theta_-,\varphi)= A_0 + \epsilon A_{in}(\varphi),\\\label{bg110}
&&\omega_{\theta}(\theta_-,\varphi)= \epsilon K_{in}(\varphi),\\\label{bg12}
&&U_2(\theta_-,\varphi)-\mu_- U_1(\th_-,\var)= \overline{U}_{2}(\theta_-)-\mu_-\overline{U}_{1}(\th_-)+ \epsilon w_1(\varphi),\\\label{bg13}
&&U_2(\theta_+,\varphi)-\mu_+ U_1(\th_+,\var)= \overline{U}_{2}(\theta_+)-\mu_+\overline{U}_{1}(\th_+)+ \epsilon w_2(\varphi),
\end{eqnarray}
where $B_{in}, A_{in}, w_1, w_2\in C^{2,\alpha}(\mathbb{T}_{2\pi})$ and $K_{in}\in C^{1,\alpha}(\mathbb{T}_{2\pi})$ for some $\alpha\in (0,1)$ and $\mu_-, \mu_+$ are two real numbers satisfying
\be\label{mu}
\mu_->0, \ \ \ e^{-\int_{\theta_-}^{\theta_+}(\overline{k}_{2}+\frac{2 \overline{A}_{12}}{\overline{A}_{22}})(\tau)d\tau} \mu_--\mu_+<0.
\ee
The conditions \eqref{bg10}, \eqref{bg11}, and \eqref{bg110} are natural due to the fact that the Bernoulli's quantity $B$, the entropy $A$, and the vorticity in the polar direction $\omega_\theta$ satisfy the transport equations \eqref{bg201} and \eqref{bg6}. These equations ensure that these quantities evolve according to the underlying flow dynamics.
Furthermore, the boundary conditions \eqref{bg12} and \eqref{bg13} impose restrictions on the flow angles at the entrance and exit. These conditions are both physically acceptable and experimentally controllable.

Now we can prove the existence and uniqueness of smooth transonic flows with nonzero vorticity that depends only on the polar and azimuthal angles.
\begin{theorem}\label{tran-exis}
{\it There exists a small positive constant $\epsilon_0$ depending only on the self-similar background solution and the boundary datum, such that if $0<\epsilon\leq \epsilon_0$, then the system \eqref{bg2} supplemented with boundary conditions \eqref{bg10}-\eqref{bg13} has a unique smooth transonic solution $({\bf U},B, A)$ with nonzero vorticity such that
\be\label{transonic3}
\sum_{j=1}^5 \|V_j\|_{C^{2,\alpha}(\overline{\mathbb{E}})}\leq C \epsilon,
\ee
where
$
V_j= U_j- \overline{U}_{j},\ \ j=1,2,\ \ V_3= U_3, \ \   V_4= B-B_0,\ \ \ V_5=A-A_0,
$
and $C$ is a positive constant depending only on $B_{in}, A_{in}, w_1, w_2$, $K_{in}$ and the background solution constructed in Proposition \ref{transonic-bg} but not on $\eps$ and $V_j$, $j=1,2,3,4,5$. Moreover, all the sonic points form a conic surface given by $\theta=s(\varphi)\in C^{2}(\mathbb{T}_{2\pi})$. The sonic surface is closed to the background sonic conic surface $\theta=\theta_{so}$ in the sense that
\be\label{estimate2}
\|s(\varphi)-\theta_{so}\|_{C^2(\mathbb{T}_{2\pi})}\leq C\epsilon.
\ee
}\end{theorem}

Before proceeding with the proof of Theorem \ref{tran-exis}, we establish a preliminary uniqueness result for the boundary value problem associated with a second-order elliptic equation.
\begin{proposition}\label{unique-elliptic}
{\it For any constant $m_0>0$, the following elliptic boundary value problem
\be\label{bg31}\begin{cases}
\p_{\theta}^2 V_1 + \frac{\overline{k}_{33}}{\sin^2\theta} \partial_{\varphi}^2 V_1 +\frac{e_1(\theta)}{\overline{A}_{22}(\theta)} V_1 + (\overline{k}_{2}+\frac{2 \overline{A}_{12}}{\overline{A}_{22}})\partial_{\theta} V_1= 0, \ \forall (\theta,\varphi)\in \mathbb{E},\\
\partial_{\theta} V_1(\theta_-,\varphi)-\mu_- V_1(\th_-,\var)=0,\ \forall\varphi\in \mathbb{T}_{2\pi},\\
\p_{\theta} V_1(\theta_+,\varphi)-\mu_+ V_1(\th_+,\var)= 0,\ \forall\varphi\in \mathbb{T}_{2\pi}.
\end{cases}\ee
has only the zero solution, where $e_1(\theta)$, $\overline{A}_{12}(\theta)$, $\overline{A}_{22}(\theta)$, $\overline{k}_{2}(\theta)$, and $\overline{k}_{33}(\theta)$ are defined on the interval $[\theta_-,\theta_+]$ as given in equations \eqref{s101}, \eqref{ba11}, \eqref{k11}, and \eqref{k33}, respectively.
}\end{proposition}

\begin{proof}

Let $V_1(\th,\var)=e^{m(\th)}\overline{V}_1(\th,\var)$, where $m(\theta)$ is a smooth function of $\theta$ to be selected, then \eqref{bg31} can be written as
\begin{align}\no
\begin{cases}
\p_\th^2\overline{V}_1+\frac{\overline{k}_{33}}{\sin^2\theta} \partial_{\varphi}^2 \overline{V}_1 + (\overline{k}_{2}+\frac{2 \overline{A}_{12}}{\overline{A}_{22}}+2m')\partial_{\theta} \overline{V}_1\\
\quad+\big(m''+(m')^2+(\overline{k}_{2}+\frac{2 \overline{A}_{12}}{\overline{A}_{22}})m' +\frac{e_1(\theta)}{\overline{A}_{22}(\theta)}\big)\overline{V}_1=0,\\
\p_\th \overline{V}_1(\th_-,\varphi)+(m'(\th_-)-\mu_-)\overline{V}_1(\th_-,\varphi)=0,\\
\p_\th \overline{V}_1(\th_+,\varphi)+(m'(\th_+)-\mu_+)\overline{V}_1(\th_+,\varphi)=0.
\end{cases}
\end{align}
Let $m_1(\theta)=m'(\theta)$, we will show that there exists a unique smooth function $m_1(\theta)$ defined on $[\theta-,\theta_+]$ such that
\begin{align}\label{m1}
\begin{cases}
m_1'(\theta)+m_1^2(\theta)+(\overline{k}_{2}+\frac{2 \overline{A}_{12}}{\overline{A}_{22}})m_1(\theta)+\frac{e_1(\theta)}{\overline{A}_{22}(\theta)}=0, &\forall \theta\in [\theta_-,\theta_+],\\
m_1(\theta_-)=\mu_->0.
\end{cases}
\end{align}
The local existence and uniqueness of $m_1(\theta)$ can be readily derived from classical ODE theory. In order to establish the boundedness of the solution $m_1$, thereby extending it to the entire interval $[\theta_-,\theta_+]$, we introduce $m_2(\theta):=a(\theta)m_1(\theta)$, which solves the following equation:
\begin{align}\nonumber
m_2'(\theta)+\frac{(m_2(\theta))^2}{a(\theta)}=-e(\theta)<0, \quad \forall \theta\in [\theta_-,\theta_+]
\end{align}
Here, $a(\theta)=e^{\int_{\theta_-}^\theta(\overline{k}_{2}+\frac{2 \overline{A}{12}}{\overline{A}_{22}})(\tau)d\tau}>0$, and $e(\theta)=\frac{e_1(\theta)a(\theta)}{\overline{A}_{22}(\theta)}>0$.
Consequently, we have $-(m_2^{-1})'(\theta)<-(a(\theta))^{-1}<0$, which implies that $m_2^{-1}(\theta)>m_2^{-1}(\theta_-)=\mu_-^{-1}>0$ for any $\theta\in(\theta_-,\theta_+]$. As a result, $0<m_2(\theta)<m_2(\theta_-)=\mu_-$, and consequently, $m_1(\theta)$ is also bounded on $[\theta_-,\theta_+]$.
Thus, we have successfully demonstrated the existence and uniqueness of a smooth solution $m_1$ to equation \eqref{m1} on $[\theta_-,\theta_+]$. Therefore, by defining $m(\theta)= \int_{\theta_-}^{\theta} m_1(\tau) d\tau$, we can conclude that $\overline{V}_1$ satisfies
\begin{align}\nonumber\begin{cases}
\p_\th^2\overline{V}_1+\frac{\overline{k}_{33}}{\sin^2\theta} \partial_{\varphi}^2 \overline{V}_1 + (\overline{k}_{2}+\frac{2 \overline{A}_{12}}{\overline{A}_{22}}+2m')\partial_{\theta} \overline{V}_1=0,\ \forall (\theta,\varphi)\in \mathbb{E},\\
\p_\th\overline{V}_1(\th_-,\var)=0,\ \forall\varphi\in \mathbb{T}_{2\pi},\\
\p_{\theta} \overline{V}_1(\th_+,\var)+(m_2(\theta_+)a(\th_+)^{-1} -\mu_+)\overline{V}_1(\theta_+,\varphi)=0,\ \forall\varphi\in \mathbb{T}_{2\pi}.
\end{cases}
\end{align}
Thanks to \eqref{mu}, $m_2(\theta_+)a(\th_+)^{-1} -\mu_+\leq \mu_-a(\th_+)^{-1} -\mu_+<0$ and by the maximum principle and Hopf's lemma, $\overline{V}_1=V_1\equiv 0$ in $\mathbb{E}$.
\end{proof}

\begin{proof}[Proof of Theorem \ref{tran-exis}.]
Denote the solution class by
\be\no
\mathcal{X}=\{{\bf V}=(V_1,\cdots V_5)^\top\in C^{2,\alpha}(\overline{\mathbb{E}}): \sum_{j=1}^5 \|V_j\|_{C^{2,\alpha}(\overline{\mathbb{E}})}\leq \delta\},
\ee
where the parameter $\delta$ will be specified later. For any $\overline{{\bf V}}\in \mathcal{X}$, we construct a mapping $\mathcal{M}$ from $\mathcal{X}$ to itself as follows.

We first solve the Bernoulli's function and the entropy
\begin{eqnarray}\no\begin{cases}
(\overline{U}_{2}+\overline{V}_2) \partial_{\theta} V_4+ \frac{\overline{V}_3}{\sin\theta}\partial_{\varphi} V_4=0,\\
(\overline{U}_{2}+\overline{V}_2) \partial_{\theta} V_5+ \frac{\overline{V}_3}{\sin\theta}\partial_{\varphi} V_5=0,\\
V_4(\theta_-,\varphi)= \epsilon B_{in}(\varphi),\ \ V_5(\theta_-,\varphi)= \epsilon A_{in}(\varphi).
\end{cases}\end{eqnarray}
By the characteristic method, there exists a unique smooth solution $(V_4,V_5)\in C^{2,\alpha}(\overline{\mathbb{E}})$ with the estimate
\be\label{bg141}
\sum_{j=4}^5\|V_j\|_{C^{2,\alpha}(\overline{\mathbb{E}})}\leq C\epsilon(\|(B_{in},A_{in})\|_{C^{2,\alpha}(\mathbb{T}_{2\pi})})\leq C_1\epsilon.
\ee
where $C_1$ is a positive constant depending only on $B_{in}, A_{in}, w_1, w_2$, $K_{in}$ and the background solution constructed in Proposition \ref{transonic-bg} but not on $\delta$, $\eps$ and $V_j$, $j=1,2,3,4,5$, and it may vary as the proof progresses.

For the vorticity, we first solve $\omega_{\theta}$ by the following transport equation, which is a linearized version of \eqref{bg6} near the background solution,
\begin{eqnarray}\label{bg15}\begin{cases}
\p_{\theta}\om_{\theta} + \frac{\overline{V}_3}{\overline{U}_{2}+\overline{V}_2} \frac{1}{\sin\theta} \p_{\varphi} \omega_{\theta} + \left\{\frac{\overline{U}_{1}+\overline{V}_1}{\overline{U}_{2}+\overline{V}_2}+ \frac{1}{\sin\theta}\p_{\varphi}\left(\frac{\overline{V}_3}{\overline{U}_{2}+\overline{V}_2}\right)+\frac{\cos\theta}{\sin\theta}\right\}\om_{\theta}
\\
\quad=\frac{\p_{\varphi} V_4}{(\overline{U}_{2}+\overline{V}_2)\sin\theta}-\frac{B_0+V_4-\frac{1}{2}\sum_{j=1}^3(\overline{U}_{j}+\overline{V}_j)^2}{\gamma(A_0+V_5) (\overline{U}_{2}+\overline{V}_2)}\frac{\p_{\varphi} V_5}{\sin\theta},\\
\omega_{\theta}(\theta_-,\varphi)= \epsilon K_{in}(\varphi).
\end{cases}\end{eqnarray}
Then there exists a unique smooth solution $\om_{\theta}\in C^{1,\alpha}(\mathbb{E})$ to \eqref{bg15} with
\be\no
\|\omega_{\theta}\|_{C^{1,\alpha}(\overline{\mathbb{E}})}\leq C_1 (\|K_{in}\|_{C^{1,\alpha}}\epsilon+\sum_{j=4}^5\|\p_{\varphi} V_j\|_{C^{1,\alpha}(\overline{\mathbb{E}})})\leq C_1\epsilon.
\ee
Once we have obtained the solution for $\omega_\theta$, we can define $V_2$ and $V_3$ as follows, using equation \eqref{bg7}:
\begin{eqnarray}\label{bg17}\begin{cases}
V_2:=\partial_{\theta} V_1 +  \frac{\overline{V}_3}{\overline{U}_{2}+\overline{V}_2} \omega_{\theta},\\
V_3:=\frac{1}{\sin\theta}\partial_{\varphi} U_1 -\omega_{\theta}=\frac{1}{\sin\theta}\partial_{\varphi} V_1 -\omega_{\theta}.
\end{cases}\end{eqnarray}
Then, linearize \eqref{bg301} near the background solution and use \eqref{s1}, we can derive that
\begin{eqnarray}\label{bg18}
&&(c^2(\overline{\rho})-\overline{U}_{2}^2)\p_{\theta} V_2 + \frac{c^2(\overline{\rho})}{\sin\theta}\p_{\varphi} V_3 +e_1(\theta)V_1\\\no
&&\quad\quad\quad\quad+ (e_2(\theta)-\overline{U}_{1}\overline{U}_{2})V_2 - \overline{U}_{1}\overline{U}_{2} \p_{\theta}V_1= G(V_4,\overline{{\bf V}}),
\end{eqnarray}
where $e_1(\th)$ and $e_2(\th)$ were defined in \eqref{s101}, and
\be\no
&&G(V_4,\overline{{\bf V}}):=-\big((\gamma-1)V_4-\frac{\gamma-1}{2}(\overline{V}_1^2+\overline{V}_3^2)
-\frac{\gamma+1}{2}\overline{V}_2^2\big)(\overline{U}_{2}'+\p_{\theta}\overline{V}_2)\\\no
&&\quad\q-(\gamma-1)\overline{U}_{1}\overline{V}_1\p_{\theta}\overline{V}_2+(\gamma+1)\overline{U}_{2}\overline{V}_2\p_{\theta}\overline{V}_2+\overline{V}_1\overline{V}_2(\overline{U}_{1}'+\p_{\theta}\overline{V}_1)
\\\no
&&\quad\q+(\overline{U}_{2}\overline{V}_1+\overline{U}_{1}\overline{V}_2)\p_{\theta}\overline{V}_1+(\gamma-1)(\overline{U}_{1}\overline{V}_1+\overline{U}_{2}\overline{V}_2)(2\overline{V}_1+\frac{\cos\theta}{\sin\theta}\overline{V}_2)\\\no
&&\quad\q-(\gamma-1)(V_4-\frac{1}{2}\sum_{j=1}^3 \overline{V}_j^2)\big(2(\overline{U}_{1}+\overline{V}_1)+\frac{\cos\th}{\sin\th}(\overline{U}_{2}+\overline{V}_2)\big)\\\no
&&\quad\q-\big((\gamma-1)(V_4-\overline{U}_{1}\overline{V}_1-\overline{U}_{2}\overline{V}_2-\frac12\sum_{j=1}^2\overline{V}_j^2)
-\frac{\gamma+1}{2}\overline{V}_3^2\big)\frac{1}{\sin\theta}\p_{\var}\overline{V}_3\\\no
&&\quad\q+(\overline{U}_{2}+\overline{V}_2)\overline{V}_3(\p_{\theta}\overline{V}_3+\frac{1}{\sin\theta}\p_{\varphi} \overline{V}_2)+ (\overline{U}_{1}+\overline{V}_1)\overline{V}_3\frac{1}{\sin\th}\p_{\varphi}\overline{V}_1.
\ee
By substituting equation \eqref{bg17} into \eqref{bg18} and the boundary conditions \eqref{bg12}-\eqref{bg13}, we can derive the second-order linear elliptic equation that $V_1$ must satisfy, along with its corresponding boundary conditions:
\begin{eqnarray}\label{bg19}\begin{cases}
\p_{\theta}^2 V_1 + \frac{\overline{k}_{33}}{\sin^2\theta} \partial_{\varphi}^2 V_1 +\frac{e_1(\theta)}{\overline{A}_{22}(\theta)} V_1 + (\overline{k}_{2}+\frac{2 \overline{A}_{12}}{\overline{A}_{22}})\partial_{\theta} V_1= G_1(V_4,\overline{{\bf V}}),\\
\partial_{\theta} V_1(\theta_-,\varphi)-\mu_- V_1(\th_-,\var)= \epsilon w_1(\varphi)-\frac{\epsilon\overline{V}_3(\theta_1,\varphi) K_{in}(\varphi)}{\overline{U}_{2}(\theta_-)+\overline{V}_2(\theta_1,\varphi)},\ \ \varphi\in\mathbb{T}_{2\pi},\\
\p_{\theta} V_1(\theta_+,\varphi)-\mu_+ V_1(\theta_+,\varphi)= \epsilon w_2(\varphi)-\frac{\overline{V}_3\omega_\theta}{\overline{U}_2(\theta_+)+\overline{V}_2}(\theta_+,\varphi),\ \ \varphi\in\mathbb{T}_{2\pi},
\end{cases}\end{eqnarray}
where $e_1(\theta)$, $\overline{A}_{12}(\theta)$, $\overline{A}_{22}(\theta)$, $\overline{k}_{2}(\theta)$, and $\overline{k}_{33}(\theta)$ are defined on the interval $[\theta_-,\theta_+]$ as given in equations \eqref{s101}, \eqref{ba11}, \eqref{k11}, and \eqref{k33}, respectively, and
\be\no
G_1(V_4,\overline{{\bf V}}):=\frac{G(V_4,\overline{{\bf V}})}{\overline{A}_{22}}-\p_{\theta}\big(\frac{\overline{V}_3\omega_{\theta}}{\overline{U}_{2}+\overline{V}_2} \big)-\frac{e_2-\overline{U}_{1}\overline{U}_{2}}{\overline{A}_{22}}\frac{\overline{V}_3\omega_{\theta}}{\overline{U}_{2}+\overline{V}_2} +\frac{\overline{k}_{33}}{\sin\theta}\p_{\varphi}\omega_{\theta}.
\ee
The uniqueness of classical solutions to \eqref{bg19} follows from Proposition \ref{unique-elliptic}, and then the Fredholm alternative theorem for the second-order elliptic boundary value problem implies the existence and uniqueness of a smooth solution $V_1\in C^{2,\alpha}(\overline{\mathbb{E}})$ to \eqref{bg19} with the estimate
\small\begin{align}\nonumber
&\|V_1\|_{C^{2,\alpha}(\overline{\mathbb{E}})}\\\no
&\leq C_1\big(\|G_1(V_4,\overline{{\bf V}})\|_{C^{\alpha}(\overline{\mathbb{E}})}+\|\omega_\th(\theta_+,\cdot)\|_{C^{2,\alpha}(\mathbb{T}_{2\pi})} \epsilon(\|(w_1,K_{in})\|_{C^{1,\alpha}(\mathbb{T}_{2\pi})}+\|w_2\|_{C^{2,\alpha}(\mathbb{T}_{2\pi})})\big)\\\label{bg191}
&\leq C_1(\sum_{j=1}^5\|\overline{V}_j\|_{C^{2,\alpha}(\overline{\mathbb{E}})}^2+ \|V_4\|_{C^{2,\alpha}(\overline{\mathbb{E}})}+ \epsilon)\leq C_1(\delta^2+\epsilon).
\end{align}
Thus, by using \eqref{bg17}, we can derive that
\be\nonumber
\sum_{j=2}^3 \|V_j\|_{C^{1,\alpha}(\overline{\mathbb{E}})}\leq C_1(\|\p_{\theta}V_1,\frac{1}{\sin\theta}\p_{\varphi}V_1\|_{C^{1,\alpha}(\overline{\mathbb{E}})}+\|\omega_{\theta}\|_{C^{1,\alpha}(\overline{\mathbb{E}})})\leq C_1(\delta^2 + \epsilon).
\ee
To improve the regularity of $V_2,V_3$, it follows from \eqref{bg15}, \eqref{bg17} and \eqref{bg19} that
\begin{eqnarray}\no\begin{cases}
\p_{\theta} V_3+\frac{\cos\theta}{\sin\theta} V_3-\frac{1}{\sin\theta}\p_{\varphi}V_2=-\frac{1}{\sin\theta}\big(\p_{\theta}(\omega_{\theta} \sin\theta) + \p_{\varphi}(\frac{\overline{V}_3}{\overline{U}_{2}+\overline{V}_2}\omega_{\theta})\big)\\
=\frac{\p_{\varphi} V_4}{(\overline{U}_{2}+\overline{V}_2)\sin\theta}-\frac{B_0+V_4-\frac{1}{2}\sum_{j=1}^3(\overline{U}_{j}+\overline{V}_j)^2}{\gamma(A_0+V_5) (\overline{U}_{2}+\overline{V}_2)}\frac{\p_{\varphi} V_5}{\sin\theta}\\
\q\q-\big\{\frac{\overline{U}_{1}+\overline{V}_1}{\overline{U}_{2}+\overline{V}_2}+ \frac{1}{\sin\theta}\p_{\varphi}\big(\frac{\overline{V}_3}{\overline{U}_{2}+\overline{V}_2}\big)\big\}\om_{\theta}\in C^{1,\alpha}(\overline{\mathbb{E}}),\\
\p_{\theta} V_2 + \frac{\overline{k}_{33}}{\sin\theta}\p_{\varphi} V_3= G_3(V_4,\overline{{\bf V}},V_1,\om_{\theta})\in C^{1,\alpha}(\overline{\mathbb{E}}),\\
V_2(\theta_-,\varphi)= \mu_- V_1(\theta_-,\varphi) + \epsilon w_1(\varphi),\\
V_2(\theta_+,\varphi)= \mu_+ V_1(\theta_+,\varphi) + \epsilon w_2(\varphi)
\end{cases}\end{eqnarray}
where
\be\no
G_3(V_4,\overline{{\bf V}},V_1,\om_{\theta})=-\frac{e_1}{\overline{A}_{22}}V_1-(\overline{k}_{2}+\frac{2 \overline{A}_{12}}{\overline{A}_{22}})\p_{\theta} V_1+\frac{G(V_4,\overline{{\bf V}})}{\overline{A}_{22}}
-\frac{e_2-\overline{U}_{1}\overline{U}_{2}}{\overline{A}_{22}}\frac{\overline{V}_3\omega_{\theta}}{\overline{U}_{2}+\overline{V}_2}.
\ee
Then $V_2$ and $V_3$ satisfy a first order elliptic system with normal boundary conditions on $\theta=\theta_{\pm}$. Thus $V_2,V_3\in C^{2,\alpha}(\ol{\mathbb{E}})$ and
\be\no
&&\sum_{j=2}^3 \|V_j\|_{C^{2,\alpha}(\overline{\mathbb{E}})}\leq C_1\bigg(\|\omega_{\theta}\|_{C^{1,\alpha}(\overline{\mathbb{E}})}+\sum_{j=4}^5 \|V_j\|_{C^{2,\alpha}(\overline{\mathbb{E}})}+ \|G_3(V_4,\overline{{\bf V}},V_1,\om_{\theta})\|_{C^{1,\alpha}(\overline{\mathbb{E}})}\\\label{bg194}
&&\quad\quad+\epsilon \sum_{j=1}^2\|w_j\|_{C^{2,\alpha}(\mathbb{T}_{2\pi})}+ \sum_{j=\pm}\|V_1(\theta_j,\cdot)\|_{C^{2,\alpha}(\mathbb{T}_{2\pi})}\bigg)\\\no
&&\leq C_1(\delta^2 +\epsilon).
\ee
By combining the estimates \eqref{bg141}, \eqref{bg191}, and \eqref{bg194}, we can derive the following result:
\be\no
\sum_{j=1}^5 \|V_j\|_{C^{2,\alpha}(\overline{\mathbb{E}})}\leq C_1(\delta^2+ \epsilon).
\ee
By choosing $\delta=\sqrt{\eps}$ and $0<\epsilon\leq \epsilon_1=\frac{1}{4 C_1^2}$, we can establish that $C_1(\delta^2+\epsilon)\leq \delta$. Consequently, we have shown that the mapping $\mathcal{M}$ maps $\mathcal{X}$ to itself.

Furthermore, it can be proven, with details omitted here for brevity, that $\mathcal{M}$ is contractive in $C^{1,\alpha}(\overline{\mathbb{E}})$ if $\epsilon$ is chosen to be sufficiently small.

Finally, we examine the location of all the sonic points, which satisfies $|{\bf M}(\theta,\varphi)|^2=1$, where ${\bf M}=(M_1,M_2, M_3)^\top:=\frac{1}{c(\rho,A)}(U_1,U_2,U_3)^\top$. It follows from \eqref{transonic3} that
\be\nonumber
\||{\bf M}|^2-|\ol{\bf M}|^2\|_{C^{2,\alpha}(\overline{\mathbb{E}})}\leq C\epsilon.
\ee
Note that
\be\no
|\ol{\bf M}(\theta_-)|^2>1, \ |\ol{\bf M}(\theta_+)|^2<1, \ \ \displaystyle \sup_{\theta\in[\theta_-,\theta_+]} \frac{d}{d \theta}|\ol{\bf M}(\theta)|^2<0,
\ee
thus for sufficiently small $\epsilon$,  $|{\bf M}(\theta_-,\varphi)|^2>1, \ |{\bf M}(\theta_+,\varphi)|^2<1$ for any $\varphi\in\mathbb{T}_{2\pi}$ and $\frac{\p}{\p \theta}|{\bf M}(\theta,\varphi)|^2<0$ for any $(\theta,\varphi)\in\mathbb{E}$. Therefore for each $\varphi\in \mathbb{T}_{2\pi}$, there exists a unique $s(\varphi)\in (\theta_-,\theta_+)$ such that $|{\bf M}(s(\varphi),\varphi)|^2=1$. In addition, by the implicit function theorem, the function $s\in C^{2}(\mathbb{T}_{2\pi})$. Furthermore, since
\be\no
||\ol{\bf M}(s(\varphi))|^2-|\ol{\bf M}(\theta_{so})|^2|&=&|\ol{\bf M}(s(\varphi))|^2-|{\bf M}(s(\varphi),\varphi)|^2|\\\no
&\leq& \||{\bf M}|^2-|\ol{\bf M}|^2\|_{C^{2,\alpha}(\overline{\mathbb{E}})} \leq C\epsilon,
\ee
one can deduce that $|s(\varphi)-\theta_{so}|\leq C\epsilon$ for any $\varphi\in \mathbb{T}_{2\pi}$. By differentiating the identity $|{\bf M}(s(\varphi),\varphi)|^2=1$ with respect to $\varphi$, we obtain:
\be\no
\begin{cases}
s'(\varphi)=&-\left(\p_\th|{\bf M}|^2 (s(\varphi),\varphi) \right)^{-1}\p_\var|{\bf M}|^2 (s(\varphi),\varphi)\\
s''(\varphi)=&-\left(\p_\th|{\bf M}|^2 (s(\varphi),\varphi) \right)^{-2}\bigg(\p^2_\var|{\bf M}|^2 (s(\varphi),\varphi)\p_\th|{\bf M}|^2 (s(\varphi),\varphi)\\
&\ \ \ \ \ \ \ \ \ \ \ \ \ \ \ \ \ \ \ \ \ \ \ \ \ \ \ \ \ \ \ \ \ \ \ \ \ \ \ \ \ \ -\p_\var|{\bf M}|^2 (s(\varphi),\varphi)\p^2_{\th\var}|{\bf M}|^2 (s(\varphi),\varphi)\bigg),
\end{cases}
\ee
and thus the estimate \eqref{estimate2} holds. The proof is now complete.

\end{proof}

{\bf Acknowledgement.} Part of this work was done when Weng visited The Institute of Mathematical Sciences of The Chinese University of Hong Kong in 2019. He is grateful to the institute for providing nice research environment. Weng is supported in part by National Natural Science Foundation of China 12071359, 12221001. Yuan thanks the Department of Mathematics of University of Macau for the financial support. The authors would like to express their sincere gratitude to Prof. Zhouping Xin for his insightful comments and constant support. 


\begin{thebibliography}{10}




\bibitem{bus31}
Busemann, A. Gasdynamik. Handbuch der Experimentalphysik, Vol. IV, Akademische Verlagsgesellschaft, Leipzig, 1931.


\bibitem{ccf16}
G. Chen, J. Chen and M. Feldman. {\it Transonic flows with shocks past curved wedges for the full Euler equations.} Discrete Contin. Dyn. Syst. 2016, 36(8): 4179-4211.

\bibitem{cf17}
G. Chen and B. Fang. {\it Stability of transonic shocks in steady supersonic flow past multidimensional wedges.} Advances in
  Mathematics \textbf{314} (2017), 493--539.

\bibitem{czz06}
G. Chen, Y. Zhang and D. Zhu. {\it Existence and stability of supersonic Euler flows past Lipschitz wedges.} Arch. Ration. Mech. Anal. 181, 261-310(2006).

\bibitem{ccx21}
G. Chen, J. Chen, W. Xiang. {\it Stability of attached transonic shocks in steady potential flow past three-dimensional wedges.} Comm. Math. Phys. 387 (2021), no.1, 111-138.

\bibitem{ckz21}
G. Chen, J. Kuang, Y. Zhang. {\it Stability of conical shocks in the three-dimensional steady supersonic isothermal flows past Lipschitz perturbed cones.} SIAM J. Math. Anal. 53 (2021), no.3, 2811-2862.


\bibitem{cl00}
S. Chen and D. Li. {\it Supersonic flow past a symmetrically curved cone.} Indiana Univ. Math. J. 49(4), 1411-1435 (2000).

\bibitem{ch01}
S. Chen. {\it Existence of stationary supersonic flows past a pointed body.} Arch. Ration. Mech. Anal. 156, 141-181 (2001).

\bibitem{cxy02}
S. Chen, Z. Xin, and H. Yin. {\it Global shock waves for the
  supersonic flow past a perturbed cone}, Communications in Mathematical
  Physics \textbf{228} (2002), no.~1, 47--84.

\bibitem{cf07}
S. Chen and B. Fang. {\it Stability of transonic shocks in supersonic flow past a wedge.} J. Differential Equations, 233 (1) (2007), 105-135.

\bibitem{cf48}
R. Courant and K.O. Friedrichs, {\it Supersonic flow and shock waves},
  Interscience Publishers,Inc., New York, 1948.

\bibitem{cy07}
D. Cui, H. Yin. {\it Global conic shock wave for the steady supersonic flow past a cone: isothermal case.} Pac. J. Math. 233, 257-289 (2007).

\bibitem{cy09}
D. Cui, H. Yin. {\it Global conic shock wave for the steady supersonic flow past a cone: polytropic case.} J. Differ. Equ. 246, 641-669 (2009).

\bibitem{el08}
V. Elling and T. P. Liu. {\it Supersonic flow onto a solid wedge.} Comm. Pure Appl. Math. 61 (2008), no. 10, 1347-1448.

\bibitem{ep16}
A. Enciso and D. Peralta-Salas. {\it Beltrami fields with a nonconstant proportionality factor are rare.} Arch. Ration. Mech. Anal. 220 (2016), 243-260.

\bibitem{fri58}
K. O. Friedrichs, {\it Symmetric positive linear differential equations}, Comm. Pure Appl. Math. \textbf{XI} (1958), 333-418.

\bibitem{go97}
P. Godin. {\it Global shock waves in some domains for the isentropic irrotational potential flow equations.} Commun. Partial Differ. Equ. 22, 1929-1997 (1997).



\bibitem{ku02}
Alexander~G. Kuz'min, {\it Boundary-Value Problems for Transonic Flow}, John
  Wiley ${\&}$ Sons, Ltd, West Sussex, 2002.

\bibitem{liy14}
J. Li, W. Ingo, and H. Yin. {\it On the global existence and stability of a multi-dimensional supersonic conic shock wave}, Comm. Math. Phys, \textbf{329} (2014), no.~2, 609--640.

\bibitem{lxy15}

L. Li, G. Xu and H. Yin. {\it On the instability problem of a 3-D transonic oblique shock wave.} Advances in Mathematics, Vol. 282, 2015, 443-515.




%
%
%
%




\bibitem{ll99}
W. C. Lien and T.P. Liu. {\it Nonlinear stability of a self-similar 3-dimensional gas flow.} Commun. Math. Phys. 204, 525-549 (1999).

\bibitem{llp22}
W.C. Lien, Y.Y. Liu and C. C. Peng. {\it Smooth transonic flows around cones.} Netw. Heterog. Media 17 (2022), no. 6, 827-845.

\bibitem{ma37}
J. W. Maccoll. {\it The conical shock wave formed by a cone moving at high speed.} Proceeding of the Royal Society (A) 159, 459-472 (1937).


\bibitem{mt87}
A. Majda and E. Thomann.{\it Multi-dimensional shock fronts for second-order wave equations.} Comm. Partial Differ. Equ. 12, 777-828 (1987).

\bibitem{maj90}
A. Majda. {\it One perspective on open problems in multi-dimensional conservation laws multi-dimensional hyperbolic problems and computation.} vol. 29, pp. 217-237. Springer, IMA (1990).

\bibitem{tm33}
G. I. Taylor, and J. W. Maccoll. {\it The air pressure on a cone moving at high speeds.} Proceedings of the Royal Society (A), 139,278-311 (1933).



\bibitem{wx19}
C. Wang and Z. Xin. {\it Smooth transonic flows of Meyer type in De Laval nozzles}, Arch. Ration. Mech. Anal. \textbf{232}
  (2019), no.~3, 1597--1647.

\bibitem{wx21}
C. Wang and Z. Xin. {\it Regular subsonic-sonic flows in general nozzles},  Advances in Math. {\bf 380} (2021), 107578.


\bibitem{wz09}
Z. Wang and Y. Zhang. {\it Steady supersonic flow past a curved cone}, J. Differential Equations 247 (2009), 1817-1850.

\bibitem{wex19}
S. Weng and Z. Xin. {\it A deformation-curl decomposition for three dimensional steady Euler equations (in Chinese).}
  Sci Sin Math, 2019, 49: 307-320.

\bibitem{w19}
S. Weng, A deformation-curl-Poisson decomposition to the three dimensional steady Euler-Poisson system with applications. {\it J. Differential Equations} {\bf 267} (2019) 6574-6603.

\bibitem{wxy21a}
S. Weng, Z. Xin and H. Yuan. {\it Steady compressible radially symmetric flows with nonzero angular velocity in an annulus}, J. Differential Equations 286 (2021), 433-454.

\bibitem{wxy21b}
S. Weng, Z. Xin and H. Yuan. {\it On some smooth symmetric transonic flows with nonzero angular velocity and vorticity.}, Math Models and Methods in Applied Science. Vol. 31, No. 13, pp. 2773-2817 (2021).


\bibitem{wx23a}
S. Weng, Z. Xin.  {\it Existence and stability of cylindrical transonic shock solutions under three dimensional perturbations}. arXiv:2304.02429.

\bibitem{wx23b}
S. Weng, Z. Xin. {\it Smooth transonic flows with nonzero vorticity to a quasi two dimensional steady Euler flow model.}  Arch. Ration. Mech. Anal. 248 (2024), no. 3, Paper No. 49, 62 pp.

\bibitem{wx24}
S. Weng, Z. Xin. {\it Some three dimensional smooth transonic flows for the steady Euler equations with an external force}.	arXiv:2404.18110. 

\bibitem{wz24}  
S. Weng, Y. Zhou. {\it Smooth axisymmetric transonic irrotational flows to steady Euler system with an external force in cylinders.} arXiv:2401.05029. 


\bibitem{xy06}
Z. Xin and  H. Yin. {\it Global multi-dimensional shock wave for the steady supersonic flow past a three-dimensional curved cone.} Anal. Appl. 4, 101-132 (2006).



\bibitem{xy08}
G. Xu and H. Yin. {\it Global transonic conic shock wave for the symmetrically perturbed supersonic flow past a cone.} Journal of Differential Equations \textbf{245} (2008), no.~11, 3389--3432.

\bibitem{xy09}
G. Xu and H. Yin. {\it Global multidimensional transonic conic shock wave for the perturbed supersonic flow past a cone.} SIAM J. Math. Anal. 41, 178-218 (2009).

\bibitem{xy10}
G. Xu and H. Yin. {\it Instability of one global transonic shock wave for the steady supersonic Euler flow past a sharp cone.} Nagoya J. Math. 199, 151-181 (2010).

%

\bibitem{yin06}
H. Yin. {\it Global existence of a shock for the supersonic flow past a curved wedge.} Acta Math. Sin. (Engl. Ser.) 22, 1425-1432 (2006).

\end{thebibliography}
\end{document}